\documentclass[reqno,10pt]{amsart}
\usepackage{xcolor}

\newcommand{\blue}{}

\usepackage{amssymb}
\usepackage{comment}
\usepackage{bbm}
\usepackage{esvect}
\usepackage{tikz}
\usepackage{float}
 \usepackage{graphicx}
\usepackage{fancyhdr}
\newcounter{mt}
\newcounter{mc}
\newcounter{mp}

\usepackage[all,cmtip]{xy}
\newtheorem{MainTheorem}[mt]{Theorem}
\newtheorem{MainProposition}[mp]{Proposition}
\newtheorem{Proposition}{Proposition}[section]
\newtheorem{MainCorollary}[mc]{Corollary}
\newtheorem{Definition}[Proposition]{Definition}
\newtheorem{Lemma}[Proposition]{Lemma}
\newtheorem{Theorem}[Proposition]{Theorem}
\newtheorem{Corollary}[Proposition]{Corollary}
\newtheorem{Remark}[Proposition]{Remark}
\newtheorem{Conjecture}[Proposition]{Conjecture}

\DeclareMathOperator{\inter}{int}
\renewcommand{\Re}{\mathrm {Re}}

\DeclareMathOperator{\vol}{vol}

\DeclareMathOperator{\Dens}{Dens}
\DeclareMathOperator{\Span}{Span}

\DeclareMathOperator{\Image}{Im}

\DeclareMathOperator{\Stab}{Stab}

\DeclareMathOperator{\Ker}{Ker}

\DeclareMathOperator{\GL}{GL}
\DeclareMathOperator{\PGL}{PGL}
\DeclareMathOperator{\SL}{SL}

\DeclareMathOperator{\SO}{SO}
\DeclareMathOperator{\sign}{sign}

\DeclareMathOperator{\Funk}{Funk}
\DeclareMathOperator{\RFunk}{RFunk}
\newcommand{\R}{\mathbb{R}}
\newcommand{\C}{\mathbb{C}}
\newcommand{\PP}{\mathbb{P}}

\newtheorem{exmp}{Example}[section]
\title{A Funk perspective on billiards, projective geometry and Mahler volume}
\author{Dmitry Faifman}
	\email{faifmand@tauex.tau.ac.il} 
 	\address{School of Mathematical Sciences, Tel Aviv University, Tel Aviv 6997801, Israel}
\date{\today}

\begin{document}

	\date{\today}

	\begin{abstract}
		We explore connections furnished by the Funk metric, a relative of the Hilbert metric, between projective geometry, billiards, convex geometry and affine inequalities. We first show that many metric invariants of the Funk metric are invariant under projective maps, as well as under projective duality. These include the Holmes-Thompson volume and surface area of convex subsets, and the length spectrum of their boundary, extending results of Holmes-Thompson and \'Alvarez Paiva on Sch\"affer's dual girth conjecture. We explore in particular Funk billiards, which generalize hyperbolic billiards in the same way that Minkowski billiards generalize Euclidean ones, and extend a result of Gutkin-Tabachnikov on the duality of Minkowski billiards.
		
		We next consider the volume of outward balls in Funk geometry. We conjecture a general affine inequality corresponding to the volume maximizers, which includes the Blaschke-Santal\'o and centro-affine isoperimetric inequalities as limit cases, and prove it for unconditional bodies, yielding a new proof of the volume entropy conjecture for the Hilbert metric for unconditional bodies.  As a by-product, we obtain generalizations to higher moments of inequalities of Ball and Huang-Li, which in turn strengthen the Blaschke-Santal\'o inequality for unconditional bodies. Lastly, we introduce a regularization of the total volume of a smooth strictly convex 2-dimensional set equipped with the Funk metric, resembling the O'Hara M\"obius energy of a knot, and show that it is a projective invariant of the convex body. 
%		
%		As a corollary, we obtain an elementary proof of the recently resolved volume entropy growth conjecture for the Hilbert metric.
	\end{abstract}
	
	\thanks{
	\\\indent The research was partially supported by ISF Grant 1750/20.
}
	
	\maketitle
		\tableofcontents
	%-------------------------------------------------------------------------
	
	\section{Introduction}
	
	\subsection{Funk and Hilbert metrics}
	
	The Hilbert metric in the interior of a convex set $K\subset\mathbb R^n$ is given by $$d^H_K(x,y)=\frac12 \log \frac{|ay||xb|}{|ax||yb|}, $$ where $a,x,y,b$ appear in that order on the line through $x,y$ with $a,b\in\partial K$. It is evidently invariant under invertible projective maps. Hilbert discovered this metric while attempting to generalize the Beltrami-Klein model of hyperbolic geometry. It is moreover a \emph{projective} metric, meaning that straight segments are geodesic. Hilbert's fourth problem \cite{hilbert_problems} asks for a construction of all projective metrics.
	
	A different example of a projective metric, which we now recall, was discovered by Funk \cite{funk}. It is given by the non-symmetric distance function $$d^F_K(x,y)=\log\frac{|xb|}{|yb|},$$ where $b$ is the intersection of the ray $\vv{xy}$ with $\partial K$. 
	
	Various problems in Hilbert geometry have been raised and studied, notably those of volume growth \cite{colbois_verovic, vernicos_polytopes, berck_bernig_vernicos, tholozan} and Gromov hyperbolicity \cite{benoist}. The Funk metric, on the other hand, has been for the most part denied the spotlight. In this paper we attempt to argue that at least from the perspectives of convex geometry and billiard dynamics, it is the Funk metric that is most natural to study.
	
	Let $K\subset \R^n$ be a convex compact set with non-empty interior, placed in the affine space $\R^n$. The \emph{Funk metric} on $K$ is arguably the simplest (non-symmetric) metric one can define naturally in its interior. Indeed, it is the non-reversible Finsler metric whose tangent unit ball at $x\in \inter(K)$ is $B_xK=K\subset T_x\R^n=\R^n$, namely $K$ itself with the origin fixed at $x$. 
	
	The Funk and Hilbert metrics are closely related: $$d^H_K(x,y)=\frac12(d^F_K(x,y)+d^{F}_K(y,x)).$$ It is a remarkable fact that the symmetrization of the Funk metric is already projectively invariant. This hints that the Funk metric comes close to being projectively invariant itself. We shall pursue this idea along several paths.
	
		\subsection{Summary of results}

	If $K\subset\R\PP^n$ is a convex body with interior, we may fix a hyperplane at infinity $\eta$ that is disjoint from $K$, and consider the corresponding $\eta$-Funk metric on $K$, denoted $\Funk_K^\eta$. Changing $\eta$ results in a change to the Funk metric, but of a very particular type: the Finsler norm is changed by an exact 1-form. This is the core reason behind the projective invariance of many metric invariants associated to the Funk metric, as follows.	
	
	For a submanifold $X\subset \inter(K)$, let $\Funk^\eta_K(X)$ stand for the Finsler metric that is induced on $X$ from $\Funk_K^\eta$, that is it is the intrinsic metric induced on $X$ from the Funk metric.  When $\eta$ is clear from context, e.g. when $K$ belongs to an affine space, we simply write $\Funk_K(X)$.

	\begin{MainTheorem}\label{main:projective_invariants}
		 Assume $K\subset \R^n$ is $C^2$-smooth and strictly convex. Let $\Omega\subset \inter(K)$ be a smooth $k$-dimensional submanifold with boundary. The following are invariant under the action of any invertible projective map $g$ keeping $K$ compact.
	\begin{enumerate}
	{\blue 	\item The $k$-dimensional Holmes-Thompson volume of $\Funk_K(\Omega)$, denoted $\vol^{F}_K(\Omega)$.
		\item The geodesics of $\Funk_K(\Omega)$, as well as the length of closed geodesics. 
		\item Assume $k=n$, that is $\Omega$ is a smooth domain. Then the Finsler billiard orbits in $\Funk_K(\Omega)$, as well as the length of the periodic orbits, are invariant.}
	
	\end{enumerate}
	\end{MainTheorem}

Therefore, whenever we discuss any of those invariants, we may assume $K\subset\R\PP^n$ is a convex set, without specifying an affine chart. Here and in the following, smoothness and strict convexity assumptions in statements concerning volumes can be relaxed through a continuity argument.

The \emph{reverse Funk metric} is defined similarly to the Funk metric by setting $B_xK=a_x(K)$, where $a_x$ is the antipodal map with respect to $x$, or equivalently $d_K^{RF}(x,y)=d_K^F(y,x)$. For $K, \eta$ as before we denote the corresponding reverse Funk metric on $\inter(K)$ by $\RFunk_K^\eta$. Theorem \ref{main:projective_invariants} applies equally to $\RFunk$.

In a sense, those invariants are more natural than the corresponding ones for the Hilbert metric, which is most clearly manifested by their invariance under projective duality. Recall that the \emph{polar} convex body $K^\vee\subset \PP^\vee=(\R\mathbb P^n)^\vee$ is the closure of the set of all hyperplanes disjoint from $K$.

\begin{MainTheorem}\label{main:duality} Assume $K\subset \inter(L)$ are $C^2$-smooth, strictly convex bodies in $\R\PP^n$.
	Then 	\begin{enumerate}
		\item  $\vol^F_L(K)=\vol^{F}_{K^\vee}(L^\vee)$.
		\item The billiard orbits in $\Funk_L(K)$ and $\RFunk_{K^\vee}(L^\vee)$ are in a natural bijective correspondence. Moreover, a periodic orbit corresponds to a periodic orbit, and they have equal length.
		\item The Holmes-Thompson volumes of $\Funk_L(\partial K)$ and $\Funk_{K^\vee}(\partial L^\vee)$ coincide.
		\item The geodesics of $\Funk_L(\partial K)$ and $\RFunk_{K^\vee}(\partial L^\vee)$ are in a natural bijective correspondence, and the respective closed geodesics have equal length.
	\end{enumerate}
\end{MainTheorem}

In the limit where $K$ is shrinking to a fixed point inside $L$ through homotheties, we get corresponding statements for linear spaces with Minkowski norms, as follows. 

Let $K, L\subset V=\R^n$ be fixed convex sets with $0$ in their interior, and consider $K\subset (V, \|\bullet\|_L)$ and $L^o\subset (V^*, \|\bullet\|_{K^o})$. The first statement becomes a tautology, as both volumes, when properly rescaled, approach the symplectic volume of $K\times L^o\subset V\oplus V^*$.
The second statement is the duality of the Minkowski billiards $(K, L)$ and $(L^o,-K^o)$, first established in \cite{gutkin_tabachnikov}, and which permeates the work of Artstein-Avidan and Ostrover on the relationship between Minkowski billiards, symplectic geometry and the Mahler conjecture \cite{artstein_ostrover, artstein_karasev_ostrover}, see also \cite{ostrover} for a survey. It is the approach of \cite{artstein_ostrover} that we take to establish the duality.
The third statement was established in \cite{holmes_thompson} in the particular case of $K=L$. The fourth statement is the (non-symmetric) Sch\"affer dual girth conjecture \cite{schaeffer}, which can also be understood as a duality of gliding billiard orbits. It was established in \cite{alvarez-paiva} together with the general case of the third statement.

The first statement implies an isomorphic form of duality in Hilbert geometry. Let $\vol_K^H$ be the Holmes-Thompson volume of the Hilbert metric in $K$.
\begin{MainCorollary}\label{cor:hilbert_volume} Given $K\subset \inter(L)$ in $\R\PP^n$,
	$\vol^H_L(K)\leq \frac{1}{2^n}{2n\choose n} \vol_{K^\vee}^H(L^\vee)$. 
	Equality is attained if and only if $K$ is an ellipsoid, and $L$ is a simplex.
\end{MainCorollary}
{\blue The third statement implies a duality in Hilbert geometry in the projective plane.
\begin{MainCorollary}\label{cor:hilbert_plane}
		Let $K\subset\inter(L)$ be two convex sets in $\R\PP^2$. Then the Hilbert lengths $\vol^H_L(\partial K)$ and $\vol^H_{K^\vee}(\partial L^\vee)$ coincide.
\end{MainCorollary}
}
It is not hard to verify that for an ellipsoid $B\subset\R\PP^n$, the Funk billiard  (though not the metric itself) in a domain $\Omega\subset \inter(B)$ coincides with that of the Beltrami-Klein hyperbolic metric. It follows in particular that for an ellipsoid $K\subset B\subset \R\mathbb P^n$, the Funk billiard in $K$ is completely integrable, see \cite{veselov}. We provide a Funk perspective and an alternative proof of this result in the projective plane, as follows.  

\begin{MainTheorem}\label{main:conics}
	Let $K\subset \inter(B)$ be nested ellipses in $\R\PP^2$. Then every orbit of the $B$-Funk billiard in $K$ has a caustic, which is a conic in the dual pencil defined by $\partial K, \partial B$. 
	Furthermore, if $Q$ is a caustic for $\Funk_B(K)$, and $\widehat Q\subset \mathbb P^\vee$ is the caustic of the dual orbit in $\RFunk_{K^\vee}(B^\vee)$, then $\widehat Q^*, \partial B, \partial K, Q$ form a harmonic quadruplet of conics. 
\end{MainTheorem}
Here  $Q^*\subset\mathbb P^\vee$ denotes the dual quadric, consisting of all tangent hyperplanes to $Q\subset\PP$. For a harmonic quadruplet of quadrics, see Definition \ref{def:harmonic_quadrics}. In particular, we can interpret the Poncelet porism of any pair of nested conics in terms of closed Funk (equivalently, hyperbolic) billiard orbits.

In the second part of this note, we consider the Holmes-Thompson volume of balls in the Funk metric.  Assume $K\subset\R^n$, and let $B_r^F(q, K)$ denote the outward $r$-ball in the Funk metric on $K$ centered at $q$. Denoting by $K^z\subset (\R^n)^*$ the dual body with respect to $z$, and by $\omega_n$ the volume of the Euclidean unit ball in $\R^n$, one has $\vol^{F}_K(B_r^F(0,K))=\frac{1}{\omega_n}\int_{(1-e^{-r})K}|K^z|dz$.  We show that for $K$ unconditional and any $r$, the volume of $B_r^F(0, K)$ is maximized when $K$ is an ellipsoid, which translates into the following equivalent systems of affine inequalities.

\begin{MainTheorem}\label{main:bodies_inequality}
	For an unconditional convex body $K\subset \R^n$, one has
	\begin{itemize}
	\item 	For all $0< \rho<1$
		$$ 	\int_{\rho K} |K^z|dz\leq n\omega_n^2\int_0^\rho\frac{t^{n-1}dt}{(1-t^2)^{\frac{n+1}{2}}}.$$
	\item For all $j\geq 0$,
		$$\int_{K\times K^o}\langle x, \xi\rangle^{2j}dxd\xi\leq \frac{(2\pi)^n}{n+2j}\frac{1}{ (n+2j)!!(n-2)!!(2j-1)!!} \left(\frac2\pi\right)^{\frac{1-(-1)^n}{2} }.$$
		
	\end{itemize}
	Equality in each case is attained uniquely by ellipsoids.
\end{MainTheorem}
\noindent For $j=1$, this is a theorem of Ball \cite{ball_remarks}.

In the limit $\rho\to0$, which corresponds to Funk balls of infinitesimal radius, the Blaschke-Santal\'o inequality is recovered: $$|K\times K^o|\leq |B^n\times (B^n)^o|=\omega_n^2.$$  On the other end of the scale, one has the following result of \cite{berck_bernig_vernicos}.
\begin{MainProposition}\label{prop:centro_affine_main}
	Assume $K\subset\R^n$ is $C^2$-smooth and strictly convex, and $0\in\inter(K)$. Then as $\rho\to 1^{-}$,
	$$\int_{\rho K}|K^z|dz\sim 2^{-\frac{n-1}{2}}\frac{\omega_n}{n-1}\frac{1}{(1-\rho)^{\frac{n-1}{2}}}\Omega_n(K),$$
	where $\Omega_n(K)$ is the centro-affine surface area of $K$.
\end{MainProposition}
For the definition of centro-affine surface area, see eq. \eqref{eq:centro_affine}.
Therefore in the limit of large Funk balls, Theorem \ref{main:bodies_inequality} yields the centro-affine isoperimetric inequality:
$$\Omega_n(K)\leq \Omega_n(B^n)=n\omega_n.$$
Thus Funk geometry interpolates between those two extremes through a continuum of affine inequalities. As both the Blaschke-Santal\'o and the centro-affine inequalities are valid for arbitrary convex bodies with centroid at the origin, we naturally conjecture that Theorem \ref{main:bodies_inequality} remains true in this generality. 

The Colbois-Verovic volume entropy conjecture \cite{colbois_verovic} asserts that the volume growth entropy of metric balls in Hilbert geometry is maximized by ellipsoids, or equivalently by smooth strictly convex bodies. This was recently proved by Tholozan \cite{tholozan} and Vernicos-Walsh \cite{vernicos_walsh} using very different methods, and previously established up to dimension $3$ in \cite{vernicos3d}. 
Theorem \ref{main:bodies_inequality} yields yet another proof of the volume entropy conjecture, albeit only for unconditional bodies. We write $B^H_r(q,K)$ for the metric ball in the Hilbert metric.
\begin{MainCorollary}\label{cor:entropy}
	For $K\subset \R^n$ a convex unconditional body,
		$$\limsup_{R\to \infty}\frac{1}{R}\log \vol_K^H (B^H_R(q,K))\leq n-1.$$
\end{MainCorollary}

Theorem \ref{main:bodies_inequality} follows from a more general functional form of the same inequalities. Below $\mathcal L$ is the Legendre transform, defined in eq. \eqref{eq:legendre}.
\begin{MainTheorem}\label{main:inequality_functional}
	For an unconditional Borel function $f=e^{-\phi}:\R^n\to[0,\infty)$, one has  the equivalent systems of inequalities
	\begin{itemize}
		\item For all $j\geq 0$, $$\int_{\R^n\times\R^n}\langle x,\xi\rangle ^{2j}e^{-\phi(x)-\mathcal L\phi(\xi)}dx d\xi\leq (2\pi)^n \frac{(n-2+2j)!!}{(n-2)!!(2j-1)!!}.$$
		\item For all $0\leq \rho<1$, $$\int_{\R^n\times\R^n}e^{-\phi(x)-\mathcal L \phi (\xi)+\rho\langle x,\xi\rangle}dxd\xi\leq \frac{(2\pi)^n}{(1-\rho^2)^{\frac n 2}}.$$
	\end{itemize}
	Equality in each case is attained uniquely, up to equality a.e., by $f$ that is a multiple of a gaussian.
\end{MainTheorem}
In the limit $\rho\to0$ we recover the functional Blaschke-Santal\'o inequality for unconditional functions. For $j=1$, this is a theorem of Huang and Li \cite{huang_li}. The proof we give is inspired by, and largely the same as the proof of the functional Blaschke-Santal\'o inequality of Lehec \cite{lehec}.

%Using the same methods, we also prove 
%
%\begin{MainTheorem}\label{main:inequality_functional2}
%	Let $f=e^{-\phi}:\R^n_+\to \R^n_+$ be a Borel function. Then 
%	\begin{itemize}
%		\item For all $j\geq 0$, the integral $$\int_{\R_+^n\times\R_+^n}\langle x,\xi\rangle^{j} e^{-\phi(x)-\mathcal L\phi(\xi)}dx d\xi$$ 
%		is maximized by $f(x)=\lambda e^{-\frac12 \langle Ax, x\rangle}$, where $A>0$ is diagonal.  
%		\item If $n=1$ and $0< \rho<1$,  
%		 $$\int_{\R_+\times\R_+}e^{-\phi(x)-\mathcal L\phi(\xi)}\sinh (\rho x\xi)dxd\xi\leq\frac{1}{\sqrt{1-\rho^2}}\arctan\frac{\rho}{\sqrt{1-\rho^2}},$$
%		with maximizer $f(x)=\lambda e^{-ax^2}$.
%	\end{itemize}
%	The maximizers in all cases are unique up to a.e. equality.
%\end{MainTheorem}

By Proposition \ref{prop:centro_affine_main}, the centro-affine surface area of $K$ can be viewed as a regularization of the total Holmes-Thompson volume of $K$ equipped with its Funk metric. It excels at capturing the volume growth of the metric, but has the drawback of being dependent on the precise way in which the volume is exhausted, namely through metric balls centered at a point. As a consequence, it depends on the choice of a center point; worse yet, it is not a projective invariant, defying the expectations one might entertain in light of Theorem \ref{main:projective_invariants} i). 

In an attempt to remedy this situation, in section \ref{sec:projective_total_volume} we propose a different regularization of the total volume, which turns out to be a projective invariant of the body alone. For simplicity we focus on the projective plane, although it is likely that a similar procedure can be carried out in greater generality. To state precisely, we use the standard Euclidean structure on $\R^{3}$, and locally identify $\R\PP^2$ with $S^2$. Given $\Omega\subset \inter(K)$ in $S^2$, it then holds that 
$$\vol^F_K(\Omega)=\frac{1}{\omega_n}\int_{\Omega\times K^\vee}\frac{dxd\xi}{\langle x,\xi\rangle^3},$$
where $K^\vee\subset S^2$ is the polar set given by $K^\vee=\{\xi\in S^2:\forall x\in K, \langle x,\xi\rangle\geq 0\}$.

For a convex set $K\subset S^2$ we define, borrowing a term from \cite{brylinski}, the \emph{Beta function of $K$} by the integral $$B_K(z)=\int_{K\times K^\vee}\langle x,\xi\rangle^zdxd\xi, \quad \Re z>-1.$$

\begin{MainTheorem}\label{main:beta}
	Let $K\subset S^2$ be $C^2$-smooth and strictly convex. Then $B_K(z)$ extends as a meromorphic function with simple poles contained in $\{-\frac52,-\frac72,-\frac 92,\dots\}$. Moreover, the value $B_K(-3)$ is a projective invariant of $K\subset\R\PP^2$.
\end{MainTheorem}
While the affinity to Funk volume is transparent, an explicit description of $B_K(-3)$ in terms of the Funk metric on $K$ remains to be found. The regularized total volume was inspired by the theory of M\"obius energy of knots \cite{ohara_knots}, and particularly its extensions to linking and Riesz energy \cite{ohara_solanes_linking, ohara_solanes_riesz}.

\subsection{Plan of the paper}
In section \ref{sec:preliminaries} we recall the basics of the various geometries that we use, mostly to fix notation. The rest of the paper consists of three main parts, namely sections \ref{sec:outlook}-\ref{sec:girth}, \ref{sec:every_scale}-\ref{sec:inequalities} and \ref{sec:projective_total_volume}, that are largely independent of each other. In Section \ref{sec:outlook} we dwell on the projective nature of the Funk metric, and its transformation under projective maps. Definition \ref{def:conomadic} is key. We establish the duality property of the Funk volume element, proving part i) of Theorems \ref{main:projective_invariants} and \ref{main:duality}, and Corollary \ref{cor:hilbert_volume}. In section \ref{sec:billiards} we discuss general Funk billiards and prove their projective invariance and duality properties, namely part ii) of Theorems  \ref{main:projective_invariants} and \ref{main:duality}. Then in section \ref{sec:conics} we focus on Funk billiards inside an ellipse, which is just the standard hyperbolic billiard, and prove Theorem \ref{main:conics}. In section \ref{sec:girth}, we establish the Funk analogue of Sch\"affer's dual girth conjecture, completing the proof of Theorems  \ref{main:projective_invariants} and \ref{main:duality}. In section \ref{sec:every_scale} we begin the study of the volume of balls in Funk geometry, establishing Proposition \ref{prop:centro_affine_main}. Then in section \ref{sec:inequalities} we prove the various inequalities, namely Theorems \ref{main:bodies_inequality},\ref{main:inequality_functional} and Corollary \ref{cor:entropy}. Section \ref{sec:projective_total_volume} is dedicated to the proof of Theorem \ref{main:beta}.

\subsection{Acknowledgements}
This project was born out of many fruitful discussions with Alina Stancu, whose interest and encouragement made this work possible. I greatly benefited from discussions with, and ideas contributed by Shiri Artstein-Avidan, Yaron Ostrover, Bo'az Klartag and Gil Solanes, to whom much gratitude and appreciation are extended. Thanks are also due to Constantin Vernicos for several useful comments on a first draft of the paper. The project got started during the author's stay in Montreal as a CRM-ISM postdoctoral fellow; the support provided by those institutes, as well as the excellent working atmospheres of UdeM, Concordia and McGill Universities, are gratefully acknowledged.

\section{Preliminaries}\label{sec:preliminaries}

\subsection{Convexity}
A convex body $K\subset V=\R^n$ is a compact convex set, which we will henceforth assume to have non-empty interior. We write $|K|$ for the Euclidean volume of $K$, and $\mathcal H^{n-1}$ for the Hausdorff measure on $\partial K$.

The Euclidean unit ball is $B^n$, and $\omega_n:=|B^n|=\frac{\pi^{\frac n 2}}{\Gamma(\frac n 2 +1)}$. The support function of $K\subset V$ is $h_K(\xi)=\sup_{x\in K}\langle x,\xi\rangle:V^*\to \R$, which is a convex function. 
We say that $K$ is smooth and strictly convex if $\partial K$ is $C^2$-smooth of strictly positive gaussian curvature.

By $\mathcal K_0(V)$ we denote the class of convex bodies with $0\in\inter(K)$. For $K\in\mathcal K_0(V)$, the \emph{Minkowski functional} of $K$ is $\|x\|_K=\inf\{t>0: \frac{x}{t}\in K\}$. If $K=-K$, it is a norm for which $K$ is the unit ball. 
 The cone measure on $\partial K$ is denoted $\mu_K$, and has $\int_K d\mu_K=n|K|$. If $\partial K$ is $C^1$ smooth, it is given by $d\mu_K(x)=\langle x, \nu_x\rangle d\mathcal H^{n-1}(x)$, where $\nu_x$ is the outward unit normal.	

The \emph{polar} (or \emph{dual}) convex body is $K^o=\{\xi\in V^*: \langle \xi,x\rangle\leq 1\}\in\mathcal K_0(V^*)$. It holds that $\|\xi\|_{K^o}=h_K(\xi)$.
For $z\in \inter(K)$, we write $K^z=(K-z)^o$.

For a Borel function $\phi:V\to \R\cup\{+\infty\}$, its \emph{Legendre transform} is \begin{equation}\label{eq:legendre}\mathcal L\phi(\xi)=\sup_{x\in\R^n} (\langle x,\xi\rangle -\phi(x)): V^*\to \mathbb R\cup\{+\infty\}.
\end{equation}
If $\phi$ is defined on a subset of $V$, we first extend it by $+\infty$.

The function $\mathcal L\phi$ is always convex, and if $\phi$ is convex then $\mathcal L^2\phi=\phi$.
An easy computation shows that $\mathcal L(\frac1p\|x\|_K^p)=\frac1q\|\xi\|_{K^o}^q$ for  $p^{-1}+q^{-1}=1$. This includes the case $p=\infty, q=1$ which reads $\mathcal L(-\log \mathbbm 1_K)(\xi)=\|\xi\|_{K^o}=h_K(\xi)$.

Writing $\phi_z=\phi(z+\bullet)$, one immediately finds that
$$\mathcal L\phi_z(\xi)=\sup_{x} (\langle x,\xi\rangle-\phi(x+z))=\mathcal L\phi(\xi)-\langle z,\xi\rangle.$$
For much more information on convexity, we refer to \cite{schneider_book}.

\subsection{Projective geometry}
Denote $V=\R^{n+1}$. We often write $\PP=\R\PP^n=\PP(V)$ for the $n$-dimensional real projective space when the dimension is clear from context, and $\PP^\vee=\PP(V^*)$ for the dual projective space. The group of isomorphisms of the projective space is $\PGL(n+1)=\GL(n+1)/\R^*$, and its elements are the invertible projective maps, also known as homographies. 
Projective duality establishes a bijection between the points of $\mathbb P$ and the hyperplanes of $\mathbb P^\vee$, by assigning to a point $x\in\mathbb P$ the hyperplane $\mathbb P(x^\perp)\subset\PP^\vee$, also denoted $x^\perp$. This allows to identify $\PP^\vee$ with the set of hyperplanes of $\PP$, and vice versa.

A \emph{convex body} $K\subset\PP$ is the image of a closed proper convex cone in $\R^{n+1}$, namely one that does not contain a line. 
Equivalently, for any hyperplane $H\subset\PP\setminus K$, $K\subset\PP\setminus H=\R^n$ is a convex body. We will often use the same notation for both the body and the cone it defines in $V$. The \emph{polar} (or \emph{dual}) convex body is $K^\vee\subset\PP^\vee$, which is the closure of the set of all hyperplanes disjoint from $K$. $K^\vee$ is a convex body, and if $K$ is smooth and strictly convex, then so is $K^\vee$. The \emph{normal map}, also called the \emph{Legendre transform}, is the bijection $\mathcal L_K:\partial K\to\partial K^\vee$, given by $x\mapsto T_x\partial K$. The linear and projective polar bodies are closely related, see Lemma \ref{lem:two_dualities}.

Let $a, b, c, d$ lie on a line in $\PP$, in that order. Their \emph{cross ratio} is $[a, b, c, d]:=\frac{|ac||bd|}{|ab||cd|}$. It is invariant under $\PGL(n+1)$. A \emph{pencil of hyperplanes} in $\mathbb P$ is the image under projective duality of a line in $\PP^\vee$. The cross ratio of four hyperplanes on a pencil is defined as their cross ratio in $\PP^\vee$. It can be computed by intersecting the hyperplanes with a generic line, and taking the cross ratio of the respective intersection points.

A \emph{quadric} $E=[A]\subset \PP$ is given by the homogeneous quadratic equation $\{x: \langle A x, x\rangle=0\}$, for some symmetric matrix $A=A_E$. A \emph{linear pencil} of quadrics is any family of quadrics of the form $[tA_1+sA_2]$, $t, s\in\R$. 
A \emph{dual pencil} is a family of the form  $[(tA_1^{-1}+sA_2^{-1})^{-1}]$, $t, s\in\R$.

Given two non-degenerate quadrics $E_1, E_2\in\mathbb R\PP^n$, the projective map $A=E_1^{-1}E_2\in\PGL(n+1)$ is defined as follows. Let $Q_1, Q_2$ be quadratic forms on $\R^{n+1}$ such that $E_j=\{Q_j=0\}$, which are uniquely defined up to a scalar multiple. Then $A$ is represented by the endomorphism $\overline A$ given by setting $Q_2(x,y)=Q_1(\overline Ax,y)$ for all $x,y\in\R^{n+1}$.

\begin{Definition}\label{def:harmonic_quadrics}
	Four quadrics $Q_1,Q_2,Q_3,Q_4\subset\PP$ form a \emph{ nic quadruplet} if $$Q_2^{-1}Q_1=Q_4^{-1}Q_3.$$ 
\end{Definition}
We refer to e.g. \cite{projective_book} for an overview to projective geometry, and to \cite{complex_convexity} for an introduction to convexity in projective geometry.

\subsection{Finsler geometry}
A \emph{non-reversible Finsler manifold} $(M,\phi)$ is a smooth manifold equipped with a function $\phi :TM\setminus \underline 0\to [0,\infty)$, which restricts to a non-symmetric norm on each tangent space. The assumed smoothness of $\phi$ depends on the problem at hand.

The tangent unit ball at $x$ is $B_xM=\{v\in T_xM: \phi(x,v)\leq 1\}$, and the cotangent ball  is $B_x^*M=(B_xM)^o\subset T_x^*M$. The co-ball bundle is $B^*M=\cup_{x\in M}B_x^*M$.

The \emph{Holmes-Thompson} volume is the measure $\mu^{HT}:=\frac{1}{\omega_n}\pi_*(\vol_{2n}|_{B^*M})$, where $\vol_{2n}$ is the symplectic volume on $T^*M$, and $\pi:T^*M\to M$ the natural projection.

By $\Dens(V)$ we denote the one-dimensional real line of Lebesgue measures, or \emph{densities}, on $V$. A measure on a manifold can be identified with a section of the line bundle $\Dens(TM)$. The Holmes-Thompson measure has $\mu^{HT}_x(B_xM)=\frac{1}{\omega_n}|B_xM||B^*_xM|$.
For an illuminating discussion of the Holmes-Thompson volume in Finsler geometry, which plays a central role in this work, we refer to \cite{alvarez_thompson, alvarez_berck}.

Both Hilbert and Funk geometries are Finsler, and we refer to \cite{hilbert_handbook} for a comprehensive account. Let us quote a few facts we will use. 

The Funk metric of $K\subset\R^n$ is given by the Finsler norm $\phi(x,v)=\|v\|_{K-x}$. The \emph{outward ball} in the Funk geometry of $K\subset\R^n$, of radius $r$ and centered at $q$, is the set
$$B_r^F(q, K)=\{x\in \inter(K): d^F_K(q,x)\leq r.\}$$
Extrinsically, $B^F_r(q, K)=q+(1-e^{-r})(K-q)$.

Given a compact domain $\Omega\subset \inter(K)$ with piecewise smooth boundary, one may, following \cite{gutkin_tabachnikov}, consider the Funk billiard map inside $\Omega$ by requiring that whenever $x,y,z\in\partial \Omega$ are consecutive points in a billiard orbit, then $$\frac{\partial}{\partial y}(d^F(x,y)+d^F(y,z))=0.$$

\section{A projective outlook on the Funk metric} \label{sec:outlook}
\subsection{The projective co-nomadic Finsler structure}
We start with some terminology.
\begin{Definition}\label{def:conomadic}
	Two Finsler structures $F_1$, $F_2$ on $M$ are \emph{co-translate} if $F_2-F_1$ is a $1$-form on $M$.
	A \emph{co-nomadic Finsler structure} on a manifold is an equivalence class of co-translate Finsler structures.
	{\blue Similarly, we define an \emph{exact co-nomadic Finsler structure} to be the equivalence class of Finsler norms up to an exact $1$-form on $M$. }
\end{Definition}
Geometrically, a co-nomadic Finsler structure is the data of all cotangent unit balls, fixed up to translation in each cotangent space. Given a co-nomadic Finsler structure $[F]$, represented by a Finsler norm $F$, its symmetrization $F^S(v):=\frac{1}2(F(v)+F(-v))$ is clearly a  reversible Finsler metric that is independent of $F$. The Holmes-Thompson volume element of a co-nomadic Finsler structure is similarly well-defined. {\blue  Furthermore, if $N\subset M$ is a $k$-dimensional submanifold, it inherits a co-nomadic Finsler structure, in particular its $k$-dimensional Holmes-Thompson volume is well-defined.}
{\blue 
\begin{Lemma}\label{lem:exact_property}
	Given an \emph{exact} co-nomadic Finsler structure $[F]$ on a manifold $M$, represented by the Finsler norm $F$, and an oriented closed curve, its length only depends on $[F]$. Moreover, the geodesics of $F$ only depend on $[F]$.
\end{Lemma}
\proof
The first statement is clear. For the second, recall that a geodesic is any curve that locally extremizes the length functional of a curve with fixed endpoints. Modifying the norm by an exact $1$-form simply adds a constant to the length functional.
\endproof }
 Take $x\in \inter(K)$, and consider the cone $K^\vee\subset V^*$. For any non-zero $\xi\in V^*$, set $B_\xi:=K^\vee\cap (x^\perp +\xi)-\xi\subset x^\perp$. It is a convex body (projectively equivalent to $K^\vee$). Furthermore, $B_{t\xi}=tB_{\xi}$ for all $t\neq 0$. Finally if $\xi=\xi_1-\xi_2\in x^\perp$, we get $B_{\xi_2}=B_{\xi_1}+\xi$.
Consequently, $B_\bullet$ defines a convex set in $(V^*/x^\perp)^*\otimes x^\perp=T^*_x\PP$, up to translation. 
Let us call it the \emph{projective co-nomadic Finsler structure on $\inter(K)$}.

To obtain a true Finsler structure, one can proceed in several ways.
\begin{enumerate}
	\item Fix a Euclidean structure on $V$. One can then take $\xi=x^\perp$ and obtain the fixed cotangent ball $K^\vee\cap (x^\perp + x)-x\subset T^*_x\PP=T_x\PP=x^\perp$. We call it the orthogonal Finsler structure.
	\item Fix $\eta\in \inter(K^\vee)$, and take $B_x^*K:=K^\vee\cap (x^\perp +\eta_0)-\eta_0$ for any $\eta_0\in\eta\subset V^*$. As we will see, this is simply the reverse Funk metric, see Proposition \ref{prop:projective_funk}.
	\item Fix an affine-equivariant point selector $S$ in the interior of convex bodies (which need only be defined on projective images of $K^\vee$), such as the center of mass or Santal\'o point. We then get a projectively invariant construction of a Finsler metric on $\inter(K)$ which is co-translate to the reverse Funk metric.
\end{enumerate}
Given an extra input $\alpha$ from the above list, we will say that the projective co-nomadic Finsler structure is \emph{anchored} by it.

The following lemma relates the notions of linear and projective polarity. It assumes $V=\R^{n+1}$ is equipped with the standard Euclidean structure, and so all linear spaces are identified with their duals.
\begin{Lemma}\label{lem:two_dualities}
	Identify a convex body $K\subset \R\PP^n$ with the cone $\tilde K\subset V=\R^{n+1}$, and with $K_t=\tilde K\cap \{x_{n+1}=t\}$. For $x\in V\setminus \{0\}$, write $[x]=\R x\in \R\PP^n$.
	Then for $\hat p=(p,1)\in \inter(K_1)\subset\R^n$, the set $-(K_1-\hat p)^o\subset \R^n=\{x_{n+1}=0\}$ coincides with $\pi_n(K^\vee\cap ([\hat p]^\perp+e_{n+1}))$, where $\pi_n$ projects orthogonally to $\R^n$.
	\end{Lemma}
\proof
Consider $(-y,z)\in \R^n\oplus \R$ which lies in $K^\vee\cap ([\hat p]^\perp+e_{n+1})$. Then for all $(\kappa,1)\in K_1$ we have $\langle (-y, z-1), (p,1)\rangle=\langle -y,p\rangle +(z-1)=0$, and $\langle (-y,z),(\kappa,1)\rangle =\langle -y,\kappa\rangle +z\geq 0$. Thus $z=1+\langle y,p\rangle$, and $\langle y,\kappa-p\rangle \leq 1$. That is, $\pi_n(K^\vee\cap ([\hat p]^\perp+e_{n+1}))\subset \{(-y,1): y\in (K_1-\hat p)^o\}$. For the opposite inclusion, start with $(y,1)\in (K_1-\hat p)^o$, and verify that $(-y,1+\langle y, p\rangle )\in K^\vee\cap ([\hat p]^\perp+e_{n+1})$.
\endproof
In particular if $K\subset \R\PP^n$ is given in the affine chart $\{x_{n+1}=1\}=\R^n+e_{n+1}$ by $K_1\subset \R^n$, then $K^\vee\cap \{x_{n+1}=1\}=-K_1^o+e_{n+1}$.

\begin{Proposition}\label{prop:projective_funk}
	Let $K\subset \R\mathbb P^n$ be a convex set, and fix $\eta\in \inter(K^\vee)$. The $\eta$-anchored projective Finsler metric then coincides with $\RFunk_K^{\eta}$.
\end{Proposition}

\proof
Using a Euclidean structure to identify $T_x\mathbb P(V)=T_x^*\mathbb P(V)=x^\perp$, the cotangent ball of the Euclidean-anchored Finsler metric is $K^\vee\cap (x^\perp +x)$, where $x\in S^n$ is identified with $x\in \mathbb P(V)$. 

Assume $\eta=e_{n+1}$. The affine chart is $\{x_{n+1}=1\}$, where the Funk cotangent unit ball at $\hat p =(p,1)=\frac{x}{x_{n+1}}=:A(x)$ is $(K_1-p)^o\subset e_{n+1}^\perp=T_p^*K_1]$. The identification between the two tangent spaces is the differential at $x$ of the map $A:S^{n}\to \{x_{n+1}=1\}$, namely $d_xA(v)=\frac{v}{x_{n+1}}-\frac{x}{x_{n+1}^2}v_{n+1}: x^\perp \to e_{n+1}^\perp$.
The dual map is $$d_xA^*:e_{n+1}^\perp \to x^\perp,\quad  d_xA^*(u)=\frac{u}{x_{n+1}}-\frac{1}{x_{n+1}^2}\langle x,u\rangle e_{n+1}. $$
That is, $d_xA^*$ is a projection to $x^\perp$ parallel to $e_{n+1}$, followed by a $\frac{1}{x_{n+1}}$- homothety.
Thus the inverse map $(d_xA^*)^{-1}$ is just the orthogonal projection on $e_{n+1}^\perp$, followed by an $x_{n+1}$-homothety: 
$$(d_xA^*)^{-1}=x_{n+1}\pi_n.$$

Observe that $K^\vee\cap (x^\perp +e_{n+1})=x_{n+1}(K^\vee\cap (x^\perp +x))$. The image of the $\eta$-anchored cotangent ball at $x$, viewed inside $x^\perp=T_x^*\R\mathbb P^n$, is therefore

$$(x^\perp +x)\cap K^\vee -\frac{e_{n+1}}{x_{n+1}} = (x^\perp +\frac{e_{n+1}}{x_{n+1}})\cap K^\vee -\frac{e_{n+1}}{x_{n+1}}, $$
which is mapped by $(d_xA^*)^{-1}$ to 

$$ x_{n+1}\pi_n ((x^\perp +\frac{e_{n+1}}{x_{n+1}})\cap K^\vee ) =\pi _n((x^\perp+e_{n+1})\cap K^\vee)= -(K_1-\hat p)^o$$ 
using Lemma \ref{lem:two_dualities}.
\endproof
In light of the discussion following Definition \ref{def:conomadic}, Theorem \ref{main:projective_invariants} part i) follows.

\begin{Corollary}
	The symmetrization of the projective co-nomadic Finsler metric is the Hilbert metric.
\end{Corollary}
\proof
Since the Hilbert metric is the symmetrization of the Funk metric, this follows immediately from Proposition \ref{prop:projective_funk}
\endproof
\begin{Remark}
	Thus we obtain a construction of the Hilbert metric which is transparently projectively invariant and Finslerian at the same time. An equivalent description can be found in \cite[section 3]{vernicos_yang}.
\end{Remark} 

A careful examination reveals that the Funk metric itself is projectively invariant, up to the addition of an \emph{exact} 1-form.
\begin{Proposition}\label{prop:funk_exact}
	Fix $\eta, \theta\in \inter(K^\vee)$, and let $F_\eta, F_\theta$ be the corresponding Funk metrics on $\inter(K)$. Then $F_\theta-F_\eta$ is an exact 1-form on $\inter(K)$. The same holds for the reverse Funk metric.
\end{Proposition}
\proof
The statements for the Funk and reverse Funk metrics are trivially equivalent. Assume $\eta=e_{n+1}$, and use a Euclidean structure to identify $T_x^*K$ with $x^\perp$. Examining the proof of Proposition \ref{prop:projective_funk} and using the notation therein, the corresponding cotangent balls differ by a shift of $w=\frac{\theta}{\langle x,\theta\rangle}-\frac{e_{n+1}}{x_{n+1}}$. To represent this translation in the fixed affine hyperplane $\{y_{n+1}=1\}$,  we put $\hat y=(y,1)$ ,  $x=\frac{\hat y}{|\hat y|}$, $\theta=(\alpha, \beta)\in \R^n\oplus \R$. Now apply $(d_xA^*)^{-1}=x_{n+1}\pi_n$ to $w$ to get $$s(y):=(d_xA^*)^{-1}w=\frac{x_{n+1}} {\langle x,\theta\rangle}\pi_n(\theta)=\frac{1}{\langle y, \alpha\rangle +\beta}\alpha.$$
Considered as a $1$-form, $s(y)$ is exact: 
\[ s(y)=d\log(\langle y,\alpha\rangle+\beta),\]
concluding the proof. 
\endproof
{\blue  Thus the exact co-nomadic class of the Funk metric is projectively invariant.}

\subsection{The volume in Funk geometry}\label{sec:volume_funk}

For a convex set $K\subset \PP$, let us construct a smooth measure $\widetilde \mu$ on $\inter(K)$. Given $x\in \inter(K)$, choose a density $d_x=\alpha_x\otimes\nu_x\in \Dens(T_x\PP)=\Dens(x^*\otimes V/x)=\Dens(V^*/x^\perp)\otimes \Dens((x^\perp)^*)$. Then $\nu_x^*\in \Dens(x^\perp)$ is defined by having $\nu_x^*\otimes\nu_x$ the symplectic (Liouville) volume on $x^\perp\oplus (x^\perp)^*$.
Choose $\xi\in V^*/x^\perp$ with $\alpha_x(\xi)=1$, and consider the intersection in $V^*$ of the cone $K^\vee$ with the affine subspace $x^\perp+\xi$. We evaluate $c(x)=\nu_x^*(K^\vee\cap (x^\perp +\xi))$ and set $\widetilde\mu_x:=c(x)d_x\in\Dens(T_x\PP)$. It is easily verified that the choice of $\alpha_x,\nu_x$ does not matter. Thus $\widetilde\mu$ is well-defined in $\inter(K)$.

\begin{Lemma}\label{lem:holmes_thompson}
	$\frac{1}{\omega_n}\widetilde\mu$ is the Holmes-Thompson volume of the projective co-nomadic Finsler structure inside $K$.
\end{Lemma}
\proof
Immediate from the construction.
\endproof

Next we construct a measure on $\mathbb P\times \mathbb P^\vee$, which is analogous to the symplectic volume on $V\oplus V^*$. The group of projective automorphisms $\PGL(V)$ acts on $\PP^\vee$ by $g(\xi):=g^{-*}\xi$. It holds that \begin{align*}
	\Dens(T_{x,\xi}(\PP\times \PP^\vee))&=\Dens(x^*\otimes V/x\oplus \xi^*\otimes V^*/\xi)\\&=\Dens^*(x)^{n+1}\otimes\Dens^*(\xi)^{n+1}\otimes \Dens(V\oplus V^*)\\&=\Dens^*(x)^{n+1}\otimes\Dens^*(\xi)^{n+1}=\Dens^*(x\oplus \xi)^{n+1},
\end{align*}
where all equalities are equviariant for the stablizer of $(x,\xi)$ in $\PGL(V)$.

Denote the incidence manifold $\mathcal Z=\{(x,\xi): x\perp \xi \}\subset \mathbb P\times \mathbb P^\vee$. We thus arrive at
\begin{Proposition}\label{prop:invariant_measure}
	There is a one-dimensional space of projective-invariant measures on $\mathbb P\times \mathbb P^\vee\setminus \mathcal Z$, with canonic normalization. Using the standard Euclidean structure on $\R^{n+1}$ to identify $\PP$ and $\PP^\vee$ locally with the unit sphere, it is given by 
	\[d\mu(x,\xi)= |\langle x,\xi \rangle| ^{-(n+1)}d\sigma_x d\sigma_\xi\]
	where $\sigma_x,\sigma_\xi$ are the standard rotationally-invariant measures on each sphere.
\end{Proposition}

\begin{Lemma}\label{lem:two_measures}
	For a measurable set $A\subset \inter(K)$, it holds that $\widetilde\mu(A)=\mu(A\times K^\vee)$.
\end{Lemma}
\proof
We have $$\mu(A\times K^\vee)=\int_{A}dx\int_{K^\vee}\frac{d\xi}{|\langle x,\xi\rangle|^{n+1}}.$$

Next we use the Euclidean structure for all choices in the definition of $\widetilde\mu$. Assume $x=e_{n+1}$, so $T_x\mathbb P=T_x^*\mathbb P=e_{n+1}^\perp$. Choose $d_x=dx$, $\alpha_x=x$, so that $\xi=x$, and $\nu_x$, $\nu_x^*$ are both the Euclidean volume. So it remains to check that $$\int_{K^\vee\cap S^n}\frac{d\xi}{|\xi_{n+1}|^{n+1}}=|K^\vee\cap (e_{n+1}+e_{n+1}^\perp)|.$$
As the jacobian of the map \begin{equation}\label{eq:projection} S^n\to e_{n+1}+e_{n+1}^\perp,\quad \xi\mapsto \xi/\xi_{n+1}\end{equation} is $1/|\xi_{n+1}|^{n+1}$, the claim follows.
\endproof

\begin{Corollary}\label{cor:dual_volumes}
The Holmes-Thompson volume of the Funk metric in the interior of $L$, denoted $\vol^F_L$, is a projective invariant of $L$, that is independent of a choice of a hyperplane at infinity. Furthermore, if $K\subset\inter(L)$ is convex, there is a duality of Funk volumes: $\vol_L^F(K)=\vol_{K^\vee}^F(L^\vee)$.
\end{Corollary}
\proof
The first part follows from Proposition \ref{prop:projective_funk} and the paragraph after Definition \ref{def:conomadic}. The second is immediate from Lemmas \ref{lem:holmes_thompson} and \ref{lem:two_measures}.
\endproof

Let $\vol^H_K$ denote the Holmes-Thompson volume of the Hilbert metric in $\inter(K)$.
\begin{Corollary}\label{cor:hilbert_volume_sec3} Given $K\subset \inter(L)$ in $\R\PP^n$,
	$\vol^H_L(K)\leq\frac{1}{2^n}{2n\choose n} \vol_{K^\vee}^H(L^\vee)$. Equality is uniquely attained when $L$ is a simplex, and $K$ is an ellipsoid.
\end{Corollary}
\proof
Choose an affine chart and assume $0\in\inter(K)$.
By the Rogers-Shephard inequality and the duality of volumes,
\begin{align*}\vol^H_L(K)=\frac{1}{\omega_n}\int_K|\frac12 (L^z-L^z)|dz\leq \frac1{2^n\omega_n}{2n\choose n}\int_K |L^z|dz&=\frac1{2^n}{2n\choose n}\vol^F_{L}(K)\\&=\frac1{2^n}{2n\choose n}\vol^F_{K^o}(L^o).\end{align*}
On the other hand by the Brunn-Minkowski inequality,
$$\vol^F_{K^o}(L^o)=\frac{1}{\omega_n}\int_{L^o}|(K^o)^z|dz\leq \frac{1}{\omega_n}\int_{L^o}|\frac12 ((K^o)^z-(K^o)^z)|dz=\vol^H_{K^o}(L^o).$$

The Rogers-Shephard inequality becomes an equality when $L$ is a simplex, while the Brunn-Minkowski  inequality above gives equality for symmetric bodies, that is $(K^o)^z$ must have antipodal symmetry (with respect to some point) for every $z\in L^o$. 

A \emph{projective center} of a convex body in $\R\PP^n $ is a center of antipodal symmetry for the body in some affine chart containing $K$.
The sets $(K^o)^z$ are projectively equivalent to $K$, and we find that the set $L^o\subset \inter(K^\vee)$ parametrizes different choices of hyperplanes $\eta^\perp$ at infinity, for which $K\subset\R\PP^n\setminus \eta^\perp$ has antipodal symmetry with some center point. But that implies that every $\eta\in L^o$ is a projective center for $K^\vee\subset\PP^\vee$. By \cite[Theorem 9-3.]{kelly_straus}, we conclude that $K^\vee$, and therefore $K$, must be a convex quadric, that is an ellipsoid.
\endproof

\section{Funk billiards and duality}\label{sec:billiards}

\begin{figure}[h]
	\centering
	\includegraphics[width=0.4\linewidth]{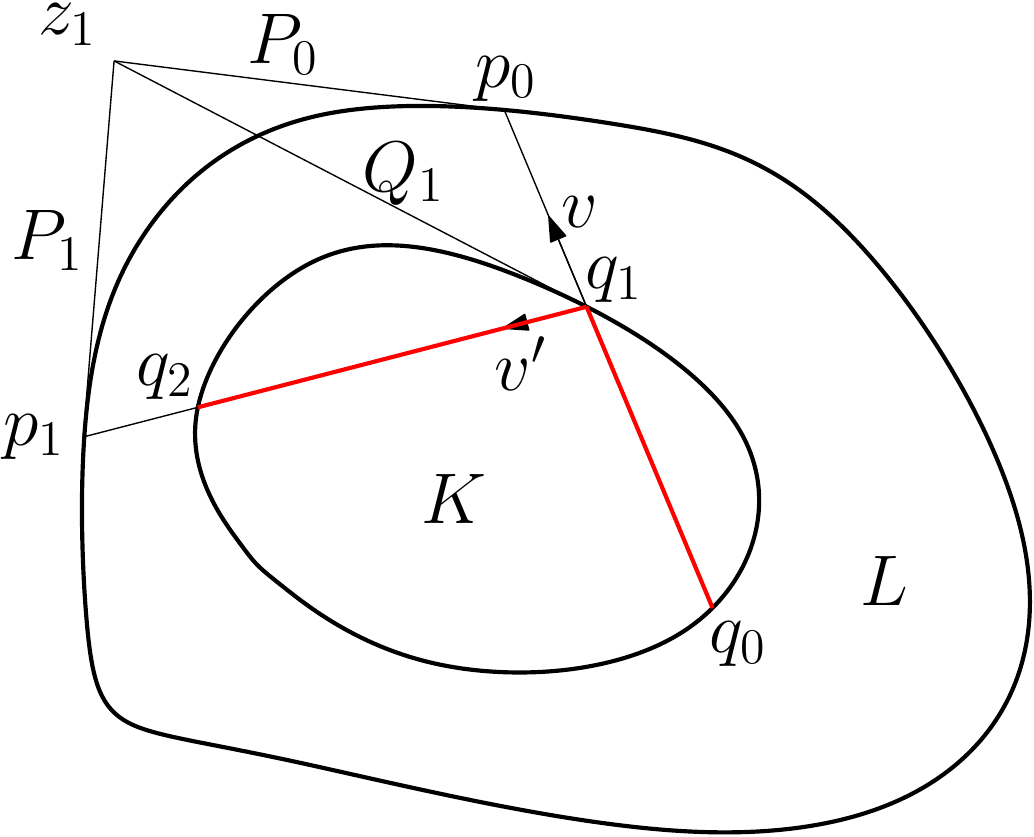}
	\caption{Funk billiard reflection law}
	\label{fig:reflectionla}
\end{figure}

Let $K\subset \inter(L)$ in $\PP^n$ be a pair of convex bodies. Using the $L$-Funk metric on $K$, we can consider the corresponding billiard dynamics inside $K$ following \cite{gutkin_tabachnikov}, which we call the $L$-Funk billiard in $K$. As it is a non-reversible Finsler structure, some care should be taken, however the results we need go through unaltered. In particular, using  \cite[Lemma 3.3]{gutkin_tabachnikov} we get the following description of the reflection law, illustrated in figure \ref{fig:reflectionla}.

Assume $q_1\in\partial K$, and $v\in S_{q_1}K$ is an incoming ray. Extend the ray until its first intersection $p_0$ with $\partial L$. Let $Q_1$ be a hyperplane tangent to $\partial K$ at $q_1$, and $P_0$ the hyperplane tangent to $L$ at $p_0$. Let $z_1=\R\PP^{n-2}$ be the intersection $Q_1\cap P_0$, and note that $z_1\subset\mathbb \R\PP^n\setminus L$. Let $P_1$ be the other tangent hyperplane to $\partial L$ through $z_1$, and $p_1\in\partial L$ the tangency point. The vector $v'\in S_{q_1}K$ pointing towards $p_1$ is then the outgoing ray.

\begin{Remark}\label{rem:proper}
It is clear from this description that a proper Funk billiard trajectory when $K, L$ are smooth and strictly convex, will extend indefinitely, in terms of bounces, as a proper billiard trajectory. 
\end{Remark}
One can also consider the billiard defined by the reverse Funk metric $\RFunk_L(K)$. Its billiard trajectories coincide with the time-reversed Funk billiard trajectories of $\Funk_L(K)$.

\begin{figure}[h]
	\centering
	\includegraphics[width=0.4\linewidth]{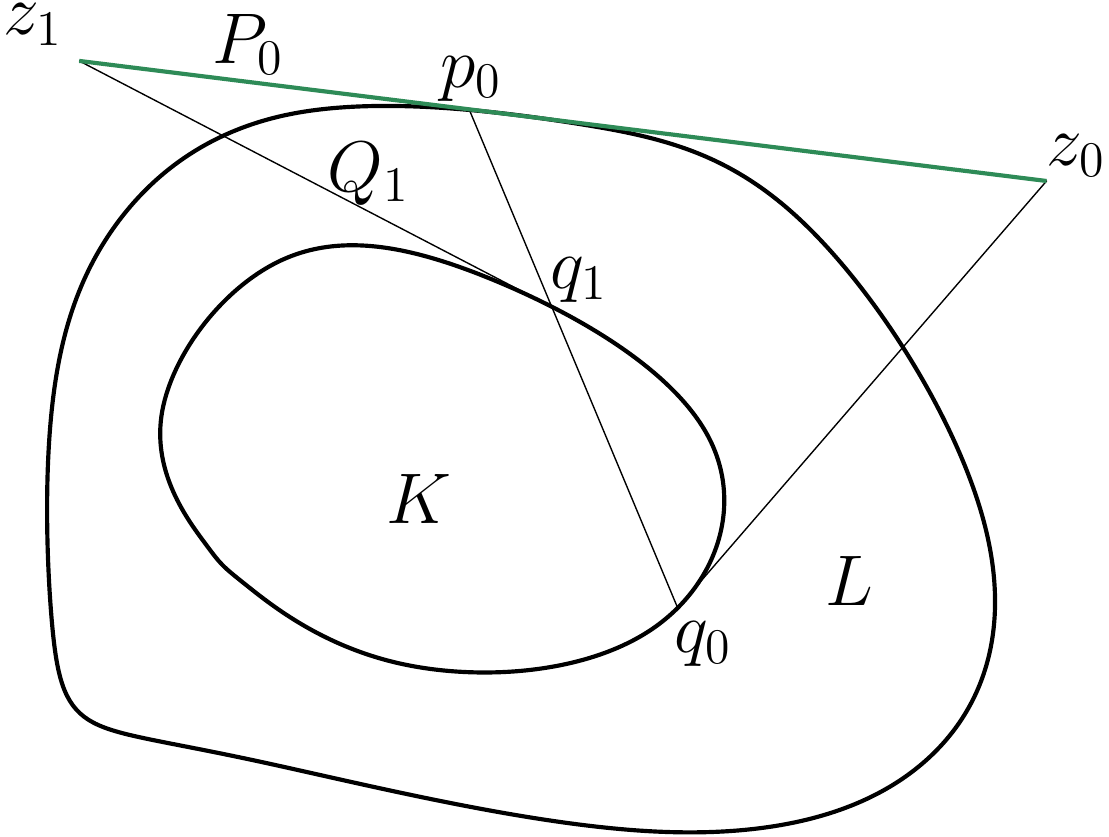}
	\caption{Funk outer billiard reflection law}
	\label{fig:outer_billiard}
\end{figure}

The same billiard can also be considered from an outer perspective: the phase space then consists of all $(n-2)$-dimensional planes $z\subset\PP^n$ that do not intersect $L$, and the billiard map mapping $z_0$ to $z_1$ is defined by the same diagram, which is reproduced in figure \ref{fig:outer_billiard} from the outer billiard perspective. By analogy with the affine setting, we call this the \emph{outer Funk billiard} on $L$ with geometry set by $K$.

\begin{figure}[h]
	\centering
	\includegraphics[width=0.4\linewidth]{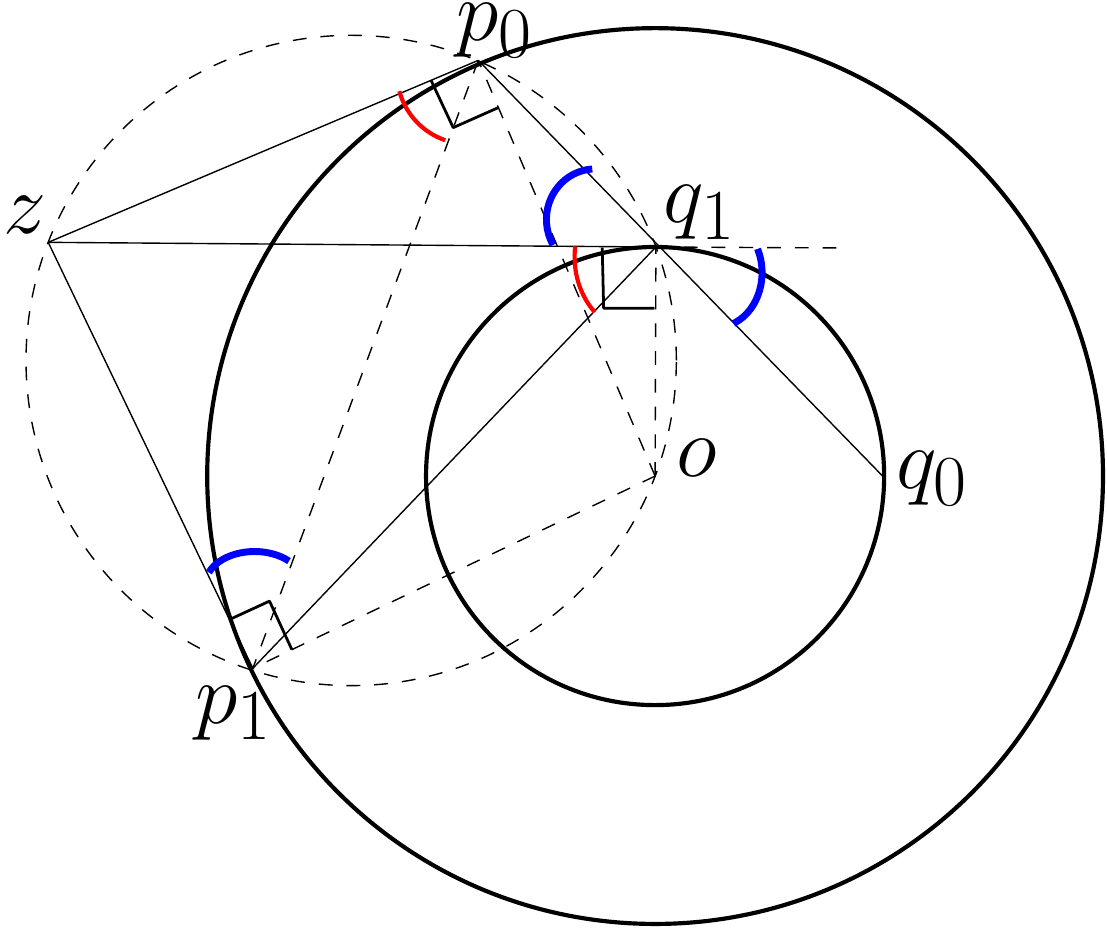}
	\caption{Funk billiard in concentric discs}
	\label{fig:circles}
\end{figure}

\begin{exmp} Assume $K\subset B\subset\R^2$ are two discs centered at the origin. Then the resulting Funk billiard inside $K$ coincides with the standard Euclidean one, as illustrated in figure \ref{fig:circles}. Indeed, take $p_0,p_1\in\partial B$, and let $z$ be the intersection of the tangents to $B$ at those points. Since $\measuredangle op_0z=\measuredangle op_1z=\measuredangle oq_1z=\frac\pi 2$, the points $p_0, p_1,q_1, z,o$ all lie on the circle having $zo$ as its diameter. It follows that the incidence angle is $\measuredangle p_0q_1z=\measuredangle p_0p_1z$, while the reflected angle is $\measuredangle p_1q_1z=\measuredangle p_1p_0z$. It remains to note that $|zp_0|=|zp_1|$.
\end{exmp}

Examining the Funk billiard law description, we see it is projectively invariant. The following provides an alternative proof of this fact, and extends the projective invariance to include the length spectrum. 

\begin{Proposition}\label{prop:billiard_invariant}
	 Choose a hyperplane at infinity $\eta\in\inter(L^\vee)$. The corresponding $L$-Funk billiard reflection law in the interior of $K$ is then independent of $\eta$, and thus projectively invariant. Furthermore, the $L$-Funk length of a periodic orbit is also projectively invariant. 
\end{Proposition}
\proof
Let $d_{\eta}(x,y)$ be the $(L,\eta)$-Funk distance from $x$ to $y$. Replacing $\eta$ by $\eta'$, we have $d_{\eta'}(x,y)=d_{\eta}(x,y)+h(y)-h(x)$ for some $h=h_{\eta,\eta'}$ by Proposition \ref{prop:funk_exact}.

Now let $p,r\in\partial K$ be fixed points. Then $p,q,r$ are consecutive reflection points of the billiard if $q$ is an extremal point for $d_\eta(p,q)+d_\eta(q,r)$. Since $d_{\eta'}(p,q)+d_{\eta'}(q,r)=d_\eta(p,q)+d_\eta(q,r)+h(r)-h(p)$, this condition is independent of $\eta$. 

Finally, if $p_0,\dots, p_{N-1}, p_N=p_0$ is a periodic orbit, then evidently 
$$\sum_{j=0}^{N-1} d_{\eta'}(p_j,p_{j+1})= \sum_{j=0}^{N-1} d_{\eta}(p_j,p_{j+1})+\sum_{j=0}^{N-1}(h(p_{j+1})-h(p_j))= \sum_{j=0}^{N-1} d_{\eta}(p_j,p_{j+1}),$$
concluding the proof.
\endproof

\begin{exmp}\label{exm:hyperbolic} Consider the Euclidean ball $B^n\subset\R^n$. The Funk Finsler norm in $\inter(B^n)$ is given by 
	$$ \phi_x(v)=\frac{\sqrt{(1-|x|^2)|v|^2+\langle x,v\rangle^2}}{1-|x|^2}+\frac{\langle x,v\rangle}{1-|x|^2}.$$
	The first summand is the Beltrami-Klein model of hyperbolic geometry . The second summand is the $1$-form $\frac{1}{1-|x|^2}\sum_{j=1}^nx_j dx_j$, which is exact. Consequently, for any domain $K\subset\inter(B^n)$, the $B^n$-Funk billiard in $K$ coincides with the hyperbolic billiard, and the length of periodic orbits is the same in both geometries.
\end{exmp}

A careful examination of the reflection law shows that it is anti-symmetric under the duality $(K,L)\leftrightarrow (L^\vee, K^\vee)$, namely a Funk billiard trajectory in one naturally defines a reverse Funk billiard trajectory in the other, which we call the dual trajectory, as follows. Let us write $Q_j=\mathcal L_K(q_j)=T_{q_j}\partial K\in\partial K^\vee$ for $q_j\in\partial K$, and similarly $P_j\in\partial  L^\vee$ for $p_j\in \partial L$. We identify the phase space with a subset of $\{(q,p): q\in \partial K, p\in \partial L^\vee\}$, where $(q,p)$ represents the outgoing ray from $q$ in the direction of $p$. 

\begin{figure}[h]
	\centering
	\includegraphics[width=1\linewidth]{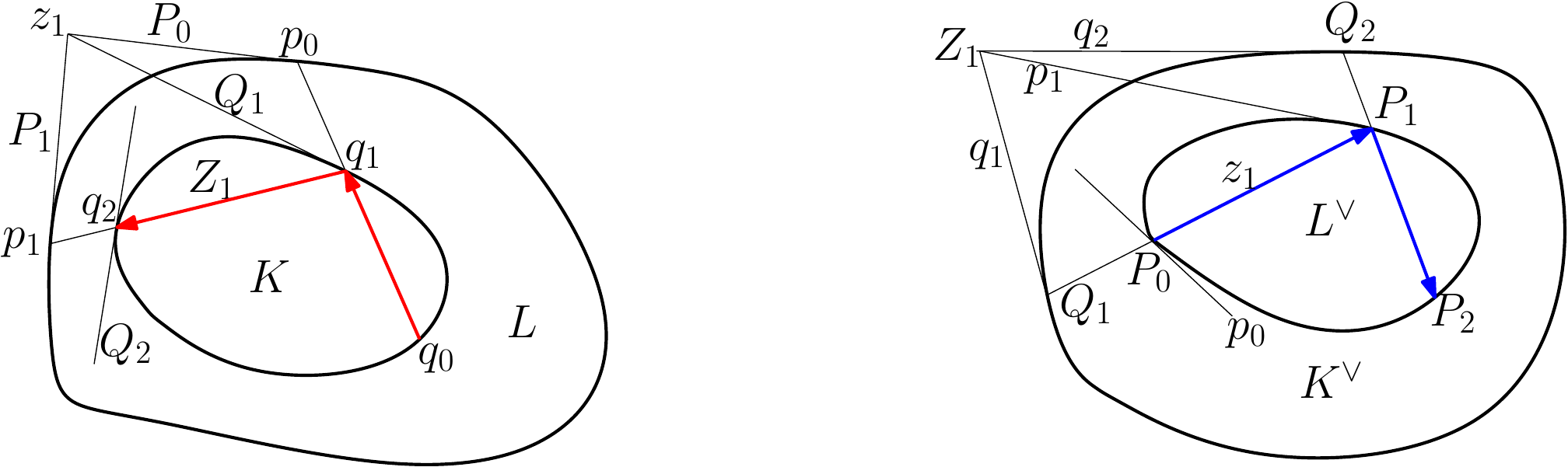}
	\caption{Dual Funk - reverse Funk billiard dynamics}
	\label{fig:dual_billiards}
\end{figure}

The billiard tranformation $(q_0,p_0)\mapsto (q_1,p_1)$ can be understood as two dual billiards unraveling simultaneously, taking alternating steps: the Funk billiard $\Funk_L(K)$, and the reverse Funk billiard $\RFunk _{K^\vee}(L^\vee)$. In $\mathbb P$, the ray arrives at $q_1\in\partial K$, and $Q_1=\mathcal L_K(q_1)$ is marked on $\partial K^\vee$. In $\mathbb P^\vee$, the interval from $Q_1$ to $P_0=\mathcal L_L(p_0)$ (which avoids $L^\vee$) is extended until its second intersection $P_1$ with $\partial L^\vee$. Now $p_1=\mathcal L_{L^\vee}(P_1)\in \partial L$ is marked, and the outgoing ray from $q_1$ is $(q_1, p_1)$. This dual dynamic is illustrated in figure \ref{fig:dual_billiards}, where the same letter is used for both a point and its dual hyperplane.

We thus established the following.

\begin{Proposition}\label{prop:duality_billiards}
	Let $q_0,q_1,\dots \in\partial K$ be an $L$-Funk billiard trajectory in $K$. Let $p_j\in\partial L$ be the first intersection of the extension of the oriented segment $[q_j, q_{j+1}]$ with $\partial L$, and $P_j$ the hyperplane tangent to $L$ at $p_j$. Then $P_0,P_1,\dots$ is a $K^\vee$-reverse Funk billiard trajectory in $L^\vee$. In the plane, the two orbits have equal rotation number.
\end{Proposition} 

\begin{Remark}
The stark similarity to the duality in Minkowski billiards is not without reason. In the limit when $L$ is much bigger than $K$, the Funk billiard becomes Minkowski, much like the hyperbolic plane is almost Euclidean in small scale. More precisely, take $K\subset L\subset\R^n$, fix $a\in L$ and set $L_R=a+R(L-a)$. Then as $R\to\infty$, the Funk billiard orbits in $K$ converge to the Minkowski billiard orbits with geometry given by $L$ centered at $a$. In particular, if $L$ is the Euclidean ball and $a$ its center, we recover the Euclidean billiard inside $K$.
\end{Remark}

Consider the case of the projective plane. Since the Funk billiard trajectory consists of straight segments and preserves the symplectic form on the phase space \cite[Theorem 4.3]{gutkin_tabachnikov}, Birkhoff's theorem on the existence of two  $m$-periodic orbits, for any $m\geq 2$ and any rotation number, remains valid.

\begin{figure}[h]
	\centering
	\includegraphics[width=0.4\linewidth]{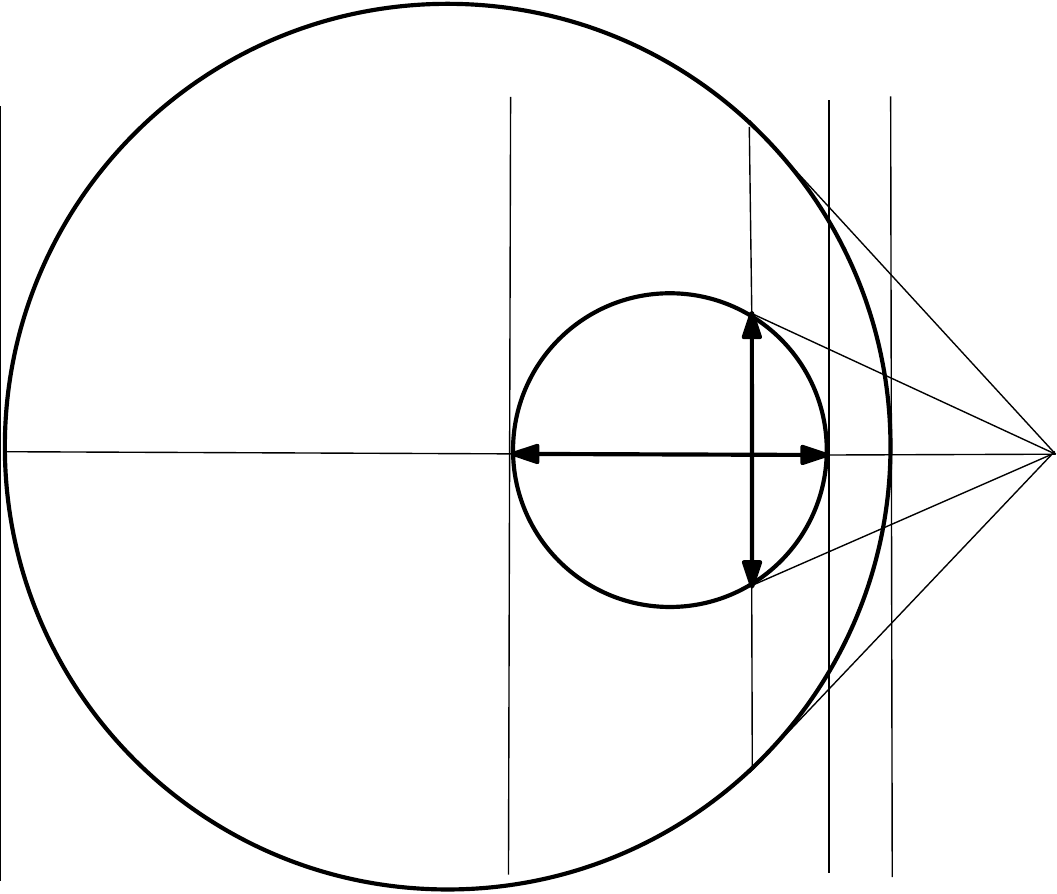}
	\caption{2-periodic orbits in a pair of discs}
	\label{fig:2_periodic_circles}
\end{figure}

As an example, the two $2$-periodic orbits when $K\subset L$ are two non-concentric circles appear in figure \ref{fig:2_periodic_circles}.

\begin{figure}[h]
	\centering
	\includegraphics[width=0.48\linewidth]{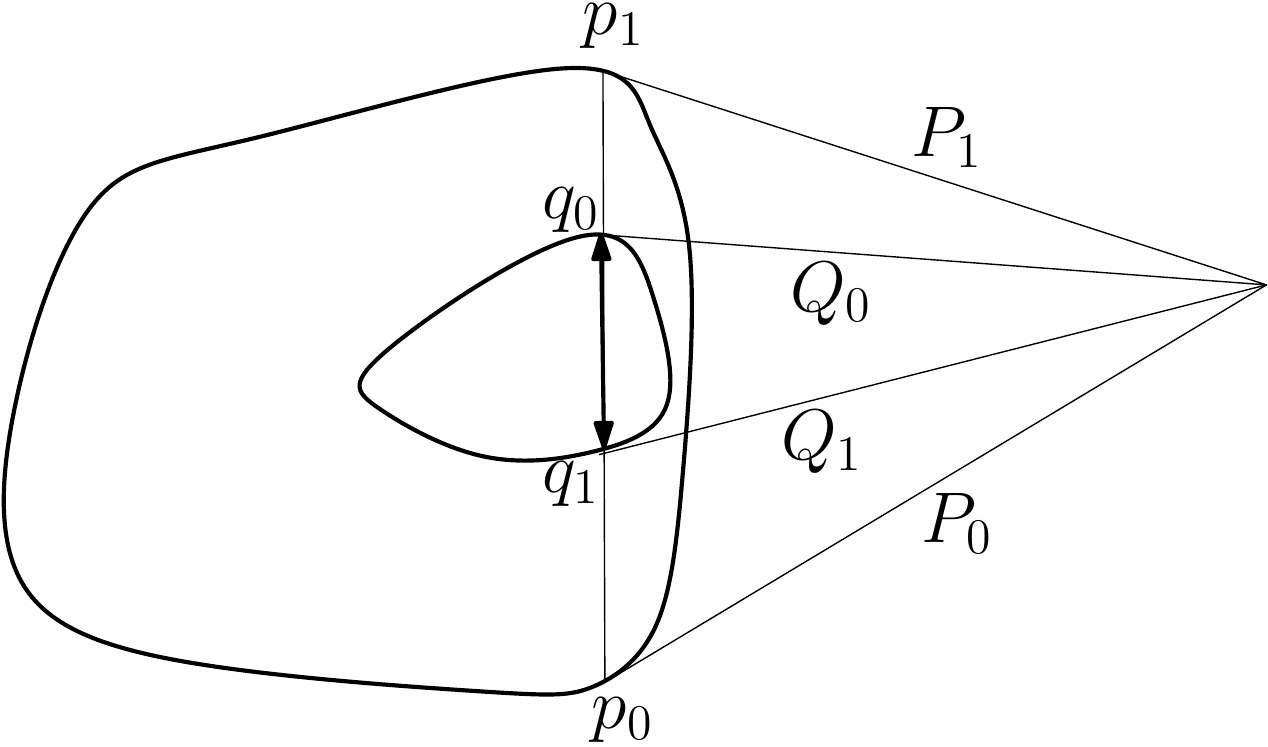}
	\caption{2-periodic orbit}
	\label{fig:2_periodic_general}
\end{figure}

To motivate the next result, let us consider 2-periodic orbits in general. Assume that $q_0, q_1\in\partial K$ are the bounce points of a 2-periodic orbit $\gamma$ of the $L$-Funk billiard in $K$, and $p_0, p_1\in \partial L$ the corresponding intersection points, appearing on the line in the order $p_1, q_0, q_1, p_0$ as in figure \ref{fig:2_periodic_general}. Note that $\gamma$ is also a reverse Funk billiard orbit. Its Funk length is twice the Hilbert distance between the points:
$$\textrm{Length}^F_L(\gamma)=\textrm{Length}^{RF}_L(\gamma)=\log \frac {|p_1q_1|}{|p_1q_0|}+\log \frac {|q_0p_0|}{|q_1p_0|}=2d^H_L(q_0, q_1).$$

Since, 
$$d^H_L(q_0,q_1)=\frac12 \log [p_1, q_0, q_1, p_0]=\frac12 \log [P_1, Q_0, Q_1, P_0]=d^H_{K^\vee}(P_0, P_1),$$
it holds that dual 2-periodic orbits have equal length. In fact, the duality of billiard orbits extends to the metric realm in full generality.

We now assume $K\subset \inter(L)$ are smooth, strictly convex bodies in $\R\PP^n$.
\begin{Theorem}\label{thm:dual_lengths}
	Dual periodic orbits in $\Funk_L(K)$ and $\RFunk_{K^\vee}(L^\vee)$ have equal lengths.
\end{Theorem}

We will need the following auxiliary fact of independent interest.

\begin{Proposition}\label{prop:unique_symplecic}
	There is a $\GL(V)$-invariant non-degenerate 2-form $\omega$ on $\mathbb P\times\mathbb P^\vee\setminus \mathcal Z$.
	Such a form is unique up to constant, and moreover it is canonically normalized. Furthermore, it is a symplectic form.
\end{Proposition}
\proof

Let $H=\Stab_{x,\xi}\subset \GL(V)$ be the stabilizer of a fixed point $(x,\xi)$, which acts on $T_x\mathbb P$ by the full linear group, and similarly on $T_\xi\mathbb P^\vee$. Observe that
$$ \wedge^2 T^*_{x,\xi}(\mathbb P\times\mathbb P^\vee)=\wedge^2(x\otimes x^\perp\oplus \xi\otimes \xi^\perp)=x^{\otimes 2}\otimes \wedge^2 x^\perp\oplus \xi^{\otimes 2}\otimes \wedge^2 \xi^\perp\oplus x\otimes x^\perp \otimes \xi\otimes\xi^\perp. $$  
The first and second summand evidently have no $H$-invariant elements. As for the third summand, we may write
$$x\otimes x^\perp \otimes \xi\otimes\xi^\perp= (x\otimes \xi)\otimes (x^\perp\otimes\xi^\perp).$$
There is a natural identification $\theta:x\otimes\xi\to \R$, and a natural pairing $\eta:x^\perp\otimes\xi^\perp \to \R$, which is moreover non-degenerate as $(x,\xi)\notin \mathcal Z$. We then let $\omega_{x,\xi}$ correspond to $\theta\otimes\eta^{-1}$.
It is immediate from the construction that $\omega$ is non-degenerate, and it is evidently canonically normalized.

For uniqueness, it remains to show that any $H$-invariant non-trivial pairing $\psi:x^\perp\otimes \xi^\perp\to \R$ must be a multiple of $\eta^{-1}$. Writing $\xi^\perp=W$ and identifying $x^\perp=W^*$, $H$ acts on $x^\perp\otimes \xi^\perp=W^*\otimes W$ by $\GL(W)$, admitting a unique invariant line. 

Finally, let us verify that $\omega$ is symplectic, that is $d\omega=0$. As $d\omega\in\Omega^3(\mathbb P\times\mathbb P^\vee\setminus \mathcal Z)$ is $\GL(V)$-invariant, it suffices to show that there are no non-trivial invariant $3$-forms.
We have an invariant decomposition
\begin{align*}
	\wedge^3 T^*_{x,\xi}(\mathbb P\times\mathbb P^\vee)&=x^{\otimes 3}\otimes \wedge^3 x^\perp\oplus \xi^{\otimes 3}\otimes \wedge^3 \xi^\perp\\&\oplus x^{\otimes 2}\otimes \wedge^2 x^\perp \otimes \xi\otimes\xi^\perp\oplus x\otimes x^\perp \otimes  \xi^{\otimes 2}\otimes \wedge^2\xi^\perp.
	\end{align*}

As before, the first and second summands do not have invariant lines. Considering the third summand, we have 
$$  x^{\otimes 2}\otimes \wedge^2 x^\perp \otimes \xi\otimes\xi^\perp=x\otimes \wedge^2 x^\perp\otimes\xi^\perp .$$
One can choose an element $g\in H$ acting by the identity on $x^\perp$ and $\xi^\perp$, and rescaling $x$ non-trivially. Thus there are no invariant elements in the third summand, and similiarly neither in the fourth.\endproof

\textit{Proof of Theorem \ref{thm:dual_lengths}.} Recall that by Proposition \ref{prop:billiard_invariant}, the Funk length of a periodic orbit is a projective invariant.
 
Now fix $\eta\in\mathbb P^\vee$, and define $H_\eta:(\mathbb P\setminus \eta^\perp)\times \mathbb P^\vee\setminus \mathcal Z\to T^*\mathbb P$ by identifying, for any $x\notin\eta^\perp$, $x^*=\eta$, and setting \begin{equation}\label{eq:identification}H_\eta(x,\xi)=-\left(\hat\eta \mapsto (x^\perp+\hat \eta)\cap \xi -\hat\eta\right) \in \eta^*\otimes x^\perp=x\otimes x^\perp=T_x^*\mathbb P.\end{equation}

For $\eta\neq \eta'$, we have in the common domain of definition the equality $$H_{\eta'}(x,\xi)=H_\eta(x,\xi)+\beta_{\eta,\eta'}(x),$$ where $\beta_{\eta,\eta'}\in\Omega^1(\mathbb P\setminus(\eta^\perp\cup\eta'^\perp))$ is a closed $1$-form by the proof of Proposition \ref{prop:funk_exact}.

Let $\omega_0$ be the canonic symplectic form on $T^*\PP$. We claim that for distinct $\eta,\eta'\in\mathbb P^\vee$, the $2$-forms $H_\eta^*\omega_0$ and $H_{\eta'}^*\omega_0$ coincide on their common domain of definition. In light of the previous paragraph, it suffices to show that $\overline\beta^*\omega_0=0$, where $\overline\beta(x,\xi)=(x,\beta(x))$, and $\beta\in \Omega^1(U)$ is a closed form, defined in a neighborhood $U\subset\mathbb P$.  We may choose coordinates locally such that $\omega_0=\sum dx_i\wedge d\xi_i$, while $\beta=\sum \frac{\partial f}{\partial x_i}dx_i$, and denote $B:=D_x\beta=(\frac{\partial^2f}{\partial x_i\partial x_j})$. Choose $u,v\in T_{x,\xi}(\mathbb P\times \mathbb P^\vee)$, and write $\hat u=dx(u), \hat v=dx(v)\in T_x\mathbb P$. We have 

\begin{align*}\overline\beta^*\omega_0(u,v)&=\omega_0((\hat u, B\hat u), (\hat v, B\hat v))\\&=\sum_i( \hat u_i (B\hat v)_i- \hat v_i (B\hat u)_i)=\sum_{i,j}B_{ij}\hat u_i\hat v_j-\sum_{i,j}B_{ji}\hat u_i\hat v_j=0.\end{align*}

It follows that as $\eta\in\mathbb P^\vee$ varies, $H_\eta^*\omega_0$ patch to a globally defined form $\omega'\in \Omega^2(\mathbb P\times\mathbb P^\vee\setminus \mathcal Z)$. Let us check that $\omega'$ is projectively invariant. First we compute
 that for $g\in \GL(V)$, \begin{align*} H_\eta\circ g(x,\xi)=H_\eta(gx,g^{-*}\xi)=(gx, \hat\eta\mapsto (g^{-*}x^\perp+\hat \eta)\cap g^{-*}\xi-\hat\eta)=gH_{g^*\eta}(x,\xi).\end{align*}
Consequently, $$g^*\omega'=g^*H_\eta^*\omega_0= H_{g^*\eta}^* g^*\omega_0= H_{g^*\eta}^*\omega_0=\omega'.$$
It follows by Proposition \ref{prop:unique_symplecic} that $\omega'=c\omega$, for some $c\in\R$. 
{\blue To find the value of $c$, we note that $\omega^n=\mu$ defines the Holmes-Thompson volume by Lemma \ref{lem:two_measures}, which by definition corresponds to $\omega_0^n$, and so $c^n=1$. It is also clear that $c>0$, e.g. by writing both forms explicitly using Euclidean coordinates at a point where $\xi=x$. Thus $c=1$.}

Replacing $T^*\mathbb P$ with $T^*\mathbb P^\vee$, one similarly has the map $$\overline{H_\theta}:\PP\times (\mathbb P^\vee\setminus \theta^\perp)\setminus \mathcal Z\to T^*\mathbb P^\vee, (x,\xi)\mapsto \left(\hat\theta \mapsto (\xi^\perp+\hat \theta)\cap x -\hat\theta\right) \in \theta^*\otimes \xi^\perp=\xi\otimes\xi^\perp, $$ and the various forms $\overline{H_\theta}^*\omega_0^\vee$ as $\theta$ varies patch to give $\omega''=\omega'=c'\omega$. The constant is $c'=1$, as the sign is flipped twice compared to $H_\eta$: once due to the change of order of $\PP, \PP^\vee$, and a second time due to the sign change in the definition of $\overline H_\theta$.

Fix $\eta\in \inter(L^\vee)$, $\theta\in \inter(K)$. It holds by Proposition \ref{prop:projective_funk} that $$S^*\Funk^\eta_L(K)=H_\eta(K\times\partial L^\vee),\quad S^*\RFunk^\theta_{K^\vee}(L^\vee)=\overline{H_\theta}(\partial K\times L^\vee).$$

Consider now $\Sigma:=\partial (K\times L^\vee)\subset\PP\times\PP^\vee\setminus\mathcal Z$. Denote by $\alpha_0, \overline\alpha_0$ the Liouville 1-forms on $T^*\mathbb P$, resp. $T^*\PP^\vee$. Define $\alpha_\eta=H_\eta^*\alpha_0$, $\overline\alpha_\theta=\overline{H_\theta}^*\overline\alpha_0$. It holds that $d\alpha_\eta-d\overline \alpha_\theta=\omega-\omega=0$, and so $\alpha_\eta-\overline \alpha_\theta=df$ for some function $f$ defined near $\Sigma$. 

A dual pair of billiard trajectories corresponds to a trajectory on $\Sigma$, as we now describe. It first follows an integral curve of $\Ker \omega|_{K\times \partial L^\vee}$ while on $\Sigma_1:=\inter(K)\times \partial L^\vee$, and of $\Ker\omega|_{\partial K\times L^\vee}$ while on $\Sigma_2 = \partial K\times \inter(L^\vee)$, corresponding to the straight segments of $\Funk_L(K)$ and $\RFunk_{K^\vee}(L^\vee)$, respectively. By Remark \ref{rem:proper}, it can only arrive to $\partial K\times\partial L^\vee$ transversally, and so the reflection law translates to simply switching from one integral curve to another, which similarly must also be transversal to $\partial K\times\partial L^\vee$. Such a trajectory on $\Sigma$ is called a generalized characteristic.

Let a generalized characteristic $\Gamma$ correspond to a dual pair of orbits $\gamma\subset \Funk_L(K)$ and $\overline\gamma\subset\RFunk_{K^\vee}(L^\vee)$. Noting that $\alpha_\eta|_{\Gamma\cap \partial K\times L^\vee}=0$, $\overline \alpha_\theta|_{\Gamma\cap K\times\partial L^\vee}=0$, we conclude that 
\[\textrm{Length}(\gamma)=\int_\Gamma \alpha_\eta=\int_\Gamma \overline\alpha_\theta=\textrm{Length}(\overline\gamma).\]\qedhere

\section{Integrability of Funk billiard in conics}\label{sec:conics}
Recall that by Example \ref{exm:hyperbolic}, the reflection laws of the $B^n$-Funk and Beltrami-Klein hyperbolic metrics coincide. Similarly, the Funk geodesics on a manifold $M\subset \inter(B^n)$ are also those of the metric induced from hyperbolic space. 

In particular, if $K$ is an ellipsoid nested in $B^n$, the Funk billiard in $K$ is completely integrable \cite{veselov, tabachnikov_hyperbolic}.
We now recover and refine this observation in the case of the projective plane by analyzing the Funk billiard directly. Let $K\subset \inter(B)$ be nested conics in $\PP^2$.

\begin{Theorem} Any given orbit in $\Funk_B(K)$ remains tangent to a conic $Q_{i}$, while the points of the outer orbit all lie on a conic $Q_o$.
	Furthermore, $Q_o$ belongs to the linear pencil of conics defined by $Q_K=\partial K,Q_B=\partial B$, while $Q_i$ belongs to the dual pencil. Moreover, the quadruplet $Q_i, Q_K, Q_B, Q_o$ is harmonic.
\end{Theorem}
\proof
We first observe that by Proposition \ref{prop:duality_billiards}, the existence of an inner caustic implies the existence of an outer caustic and vice versa. Moreover, the outer caustic belongs to the linear pencil of $Q_K, Q_B$ if and only if the inner caustic lies on the dual pencil.

Let us establish the existence of an outer caustic which belongs to the linear pencil through $Q_K$ and $Q_B$. We may assume $B$ is the unit disc. By applying a projective transformation preserving $B$, we may further assume that $Q_K$ is a centered ellipse, that is $Q_K=\{q:\langle Aq,q\rangle=1\}$. Consider a point $z_0\in B^c$, and let $z_1$ be its image under the outer Funk billiard map. Each $z_j$ belongs to a unique conic in the linear pencil, which is parametrized by $t\in \R$, $Q_t:=\{\langle (A+tI)v, v\rangle =1+t\}$. It remains to show that $t_1=t_2$, or equivalently
\begin{equation}\label{eq:t1t2} \frac{\langle Az_0,z_0\rangle -1}{1-|z_0|^2}= \frac{\langle Az_1,z_1\rangle -1}{1-|z_1|^2}.\end{equation}
Let $p_0\in Q_ B$ be the tangency point of $[z_0,z_1]$.  
Denoting $\|v\|_A:=\langle Av, v\rangle$, we have
$$ \frac{|z_0|^2-1}{|z_1|^2-1}=\frac{|z_0-p_0|^2}{|z_1-p|^2}=\frac{\|z_0-p_0\|_A^2}{\|z_1-p_0\|_A^2}, $$
the last inequality due to the fact that $z_0,pp_0,z_1$ lie on one line. Similarly, 
$$ \frac{\langle Az_0,z_0\rangle -1}{\langle Az_1,z_1\rangle -1}=\frac{\|z_0-q_0\|_A^2}{\|z_1-q_1\|_A^2}.$$
Thus \eqref{eq:t1t2} becomes 

$$ \frac{\|z_0-q_0\|_A}{\|z_0-p_0\|_A}=\frac{\|z_1-q_1\|_A}{\|z_1-p_0\|_A}\iff \frac{\sin\measuredangle_Az_0q_0p_0}{\sin\measuredangle_Az_0p_0q_0}=\frac{\sin\measuredangle_Az_1q_1p_0}{\sin\measuredangle_Az_1p_0q_1},$$
where $\measuredangle_A$ is the angle with respect to $\|\bullet\|_A$. As $\measuredangle_Az_0p_0q_0+\measuredangle_Az_1p_0q_1=\pi$, it remains to notice that $\measuredangle_Az_0q_0p_0=\frac12 \measuredangle_A q_0oq_1=\measuredangle_Az_1q_1p_0$. 

\begin{figure}[h]
	\centering
	\includegraphics[width=0.48\linewidth]{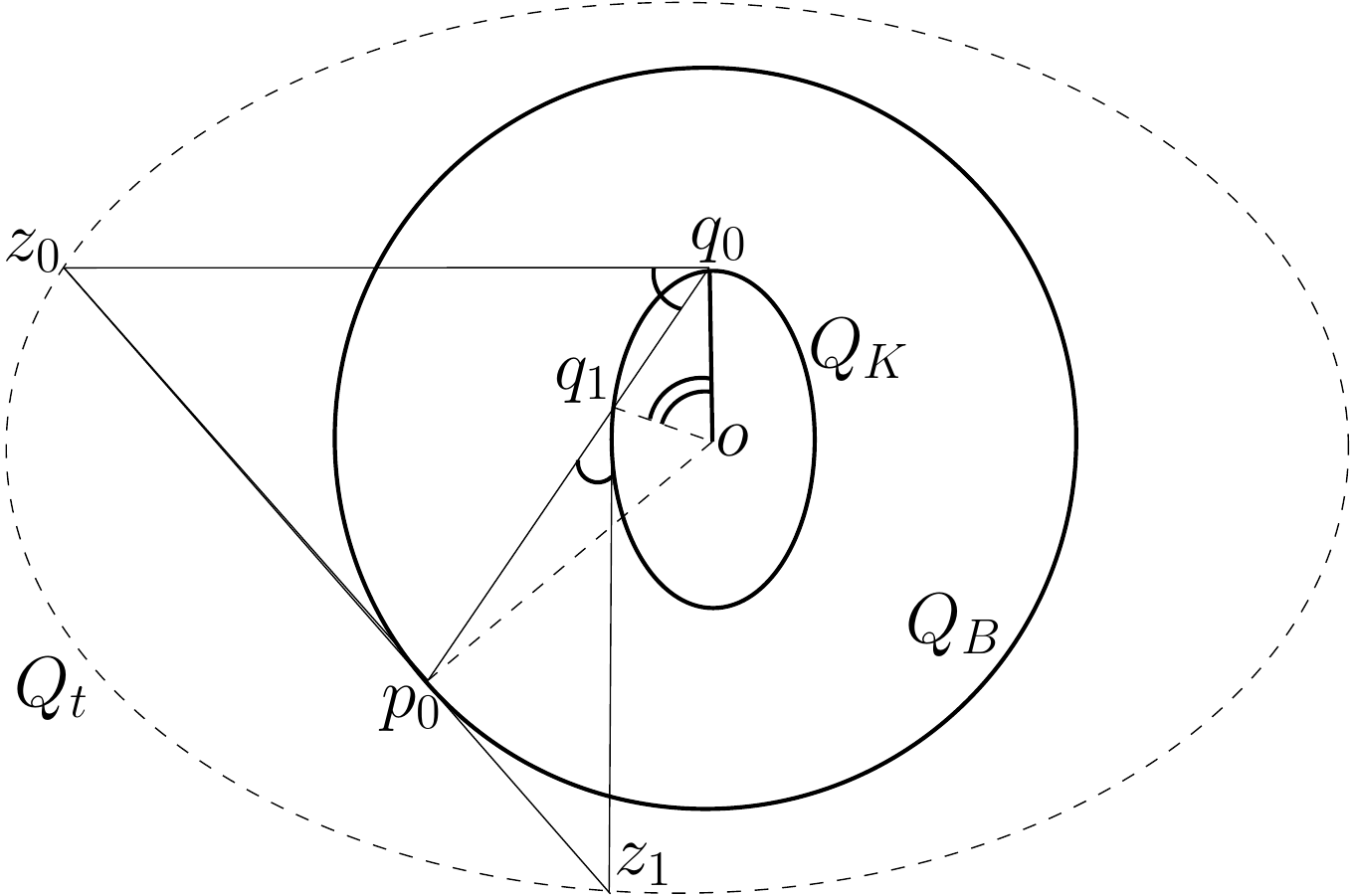}
	\caption{Outer caustic}
	\label{fig:conics1}
\end{figure}

\begin{figure}[h]
	\centering
	\includegraphics[width=0.48\linewidth]{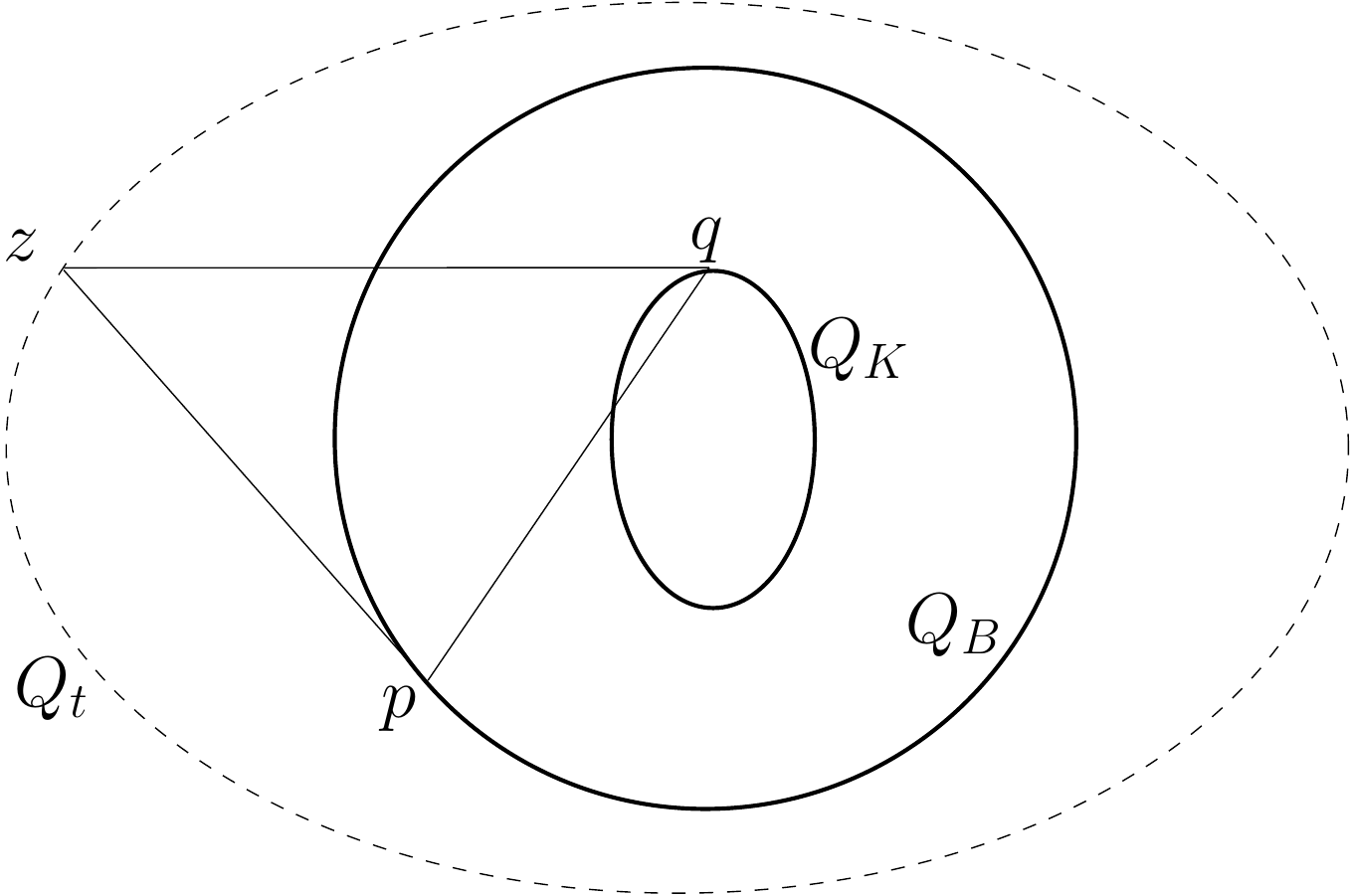}
	\caption{Relating the inner and outer caustics}
	\label{fig:conics2}
\end{figure}

For the last statement, we consider a point $z\in Q_t$, tangent to $Q_B$ at $p$ and to $Q_K$ at $q$ as in figure \ref{fig:conics2}. 

Define the conic $Q$ by $Q^{-1}:=(Q_B^{-1}Q_o)Q_K^{-1}$. To show that $Q_i=Q$, we ought to prove that $p-q$ is tangent to $Q$. Now the line through $p,q$ is represented in $\R^3$ by $w:=(p,1)\times (p-q,0)=(J(p-q), \det(p,p-q))$, where $J$ is the counterclockwise rotation by $\frac\pi 2$. It remains to check that $$\langle Q^{-1}w, w\rangle=0\iff \langle (I+tA^{-1})J(p-q), J(p-q)\rangle -(1+t)\det(p,p-q)^2=0.$$
Putting $t=-\frac{\langle Az, z\rangle-1}{|z|^2-1}$ and $ z=p+sJp=q+rJAq$ for some $s, r\in \R$, we find
$$ |z|^2-1=s^2,\quad \langle Az,z\rangle =1+r^2\langle AJAq, JAq\rangle.$$
Since $\det(p,p-q)^2+\langle p,p-q\rangle ^2=|p-q|^2$, the claimed equality becomes 
$$ s^2\langle p,p-q\rangle^2=r^2\|JAq\|_A^2(\|J(p-q)\|_{A^{-1}}^2-\det(p,p-q)^2).$$
The left hand side can be rewrriten as
\begin{align*}\langle sp, p-q\rangle^2&=\langle sJp, J(p-q)\rangle^2=\langle z-p, J(p-q)\rangle^2=\\&\langle z-q, J(p-q)\rangle^2=r^2\langle JAq, J(p-q)\rangle^2=r^2\langle Aq, p-q\rangle^2.\end{align*}
It thus remains to check that whenever $|p|=1$ and $\langle Aq,q\rangle=1$, one has
$$ \langle Aq, p-q\rangle^2=\|JAq\|_A^2(\|J(p-q)\|_{A^{-1}}^2-\det(q,p-q)^2).$$
We can write 
\begin{align*}\|J(p-q)\|_{A^{-1}}^2-\det(q,p-q)^2&=\|J(p-q)\|_{A^{-1}}^2-\langle J(p-q), q\rangle^2\\&=\|J(p-q)\|_{A^{-1}}^2-\langle Aq, J(p-q)\rangle_{A^{-1}}^2\\&=\langle J(p-q), \nu \rangle_{A^{-1}}^2\end{align*}
where $\nu=\frac{Jq}{\|Jq\|_{A^{-1}}}$, since $\|Aq\|_{A^{-1}}^2=\langle A^{-1}Aq, Aq\rangle=1$ and $\langle \nu, Aq\rangle_{A^{-1}}=0$.

It remains to check that 

$$ \|Jq\|_{A^{-1}}\langle JAq, J(p-q)\rangle=\|JAq\|_A\langle A^{-1}Jq, J(p-q)\rangle,$$
or equivalently that $$\|Jq\|_{A^{-1}}JAq=\|JAq\|_A A^{-1}Jq.$$
If $JAq=\lambda A^{-1}Jq$ then $$\lambda =\frac{\|JAq\|_A}{\|A^{-1}Jq\|_A}=\frac{\|JAq\|_A}{\|Jq\|_{A^{-1}}},$$
so it remains to check that $JA$ and $A^{-1}J$, or equivalently $AJA$ and $J$, are proportional. 
Both are anti-symmetric $2\times 2$ matrices, and so that must be the case.

\endproof
\begin{Remark}
	When $K$ approaches $B$, $Q_i$ and $Q_o$ become polars of one another with respect to $\partial B$.
\end{Remark}

\section{Sch\"affer's dual girth conjecture}\label{sec:girth}

The Gutkin-Tabachnikov duality of Minkowski billiards has a continuous counterpart, namely Sch\"affer's dual girth conjecture \cite{schaeffer}, which is a theorem of \'Alvarez Paiva \cite{alvarez-paiva}. The setting of Funk geometry is no different, providing another collection of projective invariants enjoying duality.
{\blue 
\begin{Theorem}\label{thm:girth}
Let $L$ be a smooth, strictly convex body in $\PP^n$. 
\begin{enumerate}
	\item Let $M\subset\inter(L)$ be a compact submanifold with boundary, and consider $\Funk^\eta_L(M)$. Then its Holmes-Thompson volume, its geodesics and the length of closed geodesics are all independent of $\eta$, and thus constitute projective invariants of $(M,L)$.
	\item  If $K\subset\inter(L)$ is strictly convex, then $\Funk_L(\partial K)$ and $\RFunk_{K^\vee}(\partial L^\vee)$ have the same Holmes-Thompson volumes and length spectra.
\end{enumerate}
\end{Theorem}
\proof
The first claim follows immediately by Proposition \ref{prop:funk_exact} and Lemma \ref{lem:exact_property}, as the Funk metric induces on $M$ a projectively-invariant exact co-nomadic Finsler structure.

For the second claim, we follow closely the proof in \cite{alvarez-paiva}. Write for simplicity $\partial K$ for $\Funk^\eta_L(\partial K)$, and similarly  $\partial L^\vee$ for $\RFunk^\theta_{K^\vee}(\partial L^\vee)$. 

%We first establish projective invariance. For $n=2$, the claim follows immediately from Proposition \ref{prop:funk_exact} (note also that the Holmes-Thompson volume of $\partial K$ equals the arithmetic mean of its clockwise and counterclockwise lengths, that is its Hilbert length). 
%For $n=2$, the Holmes-Thompson volume of $\partial K$ equals the arithmetic mean of its clockwise and counterclockwise Funk lengths, that is its Hilbert length. The duality in this case can be seen as a limit of Theorem \ref{thm:dual_lengths}, where we use Birkhoff's existence theorem to consider dual pairs of $m$-periodic orbits of rotation number $1$ with $m\to\infty$.

%Thus we may assume $n\geq 3$, in particular any closed 1-form on the tangent sphere bundle of $S^{n-1}$ is exact. This can be checked separately for $n=3$, when $S(S^3)$ is diffeomorphic to $\R\PP^3$ which has $H^1(\R\PP^3, \R)=H^1(S^3, \R)=0$, and for $n\geq 4$ by the long exact sequence of homotopy groups associated to the sphere bundle: 
%$$ \pi_2(S^{n-1})\to \pi_1(S^{n-2})\to \pi_1(S(S^{n-1}))\to \pi_1(S^{n-2}),$$ 
%and $ \pi_2(S^{n-1})= \pi_1(S^{n-2})=\pi_1(S^{n-2})=0$.

Define a smooth embedding $\Phi_\eta:S^*\partial K\to \partial K\times\partial L^\vee$ as follows. Let $\pi_x: T_x^*K\to T_x^*\partial K$ be the natural projection. Recall that by Proposition \ref{prop:projective_funk}, $S^*_xK=H_\eta(x\times \partial L^\vee)$, where $H_\eta: K\times L^\vee\to T^*K$ is given by eq. \eqref{eq:identification}. Now $\pi_x(S^*_xK)=B_x^*\partial K$, and $S^*_x\partial K$ is diffeomorphic to its preimage under $\pi_x$, namely the shadow boundary $Z_x\subset S_x^* K$. Denote $\tau_x=\pi_x^{-1}:S^*_x\partial K\to Z_x$, and define  $\tau_\eta:S^*\partial K\to S^*K, (x,\xi')\mapsto (x,\tau_x(\xi'))$. Set $$\Phi_\eta(x,\xi')=H_\eta^{-1}\circ\tau_\eta(x, \xi').$$

Let us describe the image of $\Phi_\eta$. It holds that $(x,\xi)\in \Image(\Phi_\eta)$ if and only if the restricted projection $\pi_x:T_\xi S^*_xK\to T_x^*\partial K$ is not an isomorphism, or equivalently if the natural pairing $T_x\partial K\otimes T_\xi\partial L^\vee\to \R$ is degenerate. The latter pairing is induced from the natural non-degenerate pairing $T_x\mathbb P\otimes T_\xi \mathbb P^\vee\to \R$, which is given by identifying $\xi=x^*$, $(x^\perp)^*=\xi^\perp$, so that $T_x \PP=x^*\otimes (x^\perp)^*=\xi\otimes \xi^\perp=T_\xi^*\PP^\vee$. In particular, $\Image(\Phi_\eta)$ does not depend on $\eta$.

Turning to the dual manifold, define $\overline{\Phi_\theta}=\overline{H_\theta}^{-1}\circ\tau_\theta:S^*\partial L^\vee\to \partial K\times\partial L^\vee$ in full analogy with $\Phi_\eta$. Observe that the images of the two maps $\Phi_\eta$, $\overline \Phi_\theta$ coincide by their common description as the degeneracy locus of the pairing  $T_x\partial K\otimes T_\xi\partial L^\vee\to \R$.

%Recall that by the proof of Theorem \ref{thm:dual_lengths}, $H_\eta^*\omega_0=\omega$, where $\omega_0$ is the symplectic form of $T^*K$, and $\omega$ as in Proposition \ref{prop:unique_symplecic}. 
%Let $\alpha_1$ be the canonic contact form on $S^*\partial K$. We will show that $\Phi_\eta^*\omega=d\alpha_1$. 
%Since $H_\eta\circ \Phi_\eta=\tau_\eta$, we have
%$$\Phi_\eta^*\omega=\Phi_\eta^*H_\eta^*\omega_0=\tau_\eta^*\omega_0,$$
%so it remains to verify that $\tau_\eta^*\omega_0=d\alpha_1$. Now $T^*\PP$ has the canonic Liouville 1-form $\alpha_0$, and $d\alpha_0=\omega_0$, so it suffices to check that $\tau_\eta^*\alpha_0=\alpha_1$. But this is evident: $\alpha_1|_{x,\xi'}(u)=\xi'(dx(u))=\tau_x(\xi')(dx(u))=\tau_\eta^*\alpha_0(u)$, since $dx(u)\in T_x\partial K$.

Let $\alpha_1$ be the canonical contact form on $S^*\partial K$, and $\overline\alpha_1$ that of $S^*\partial L^\vee$. 
%Denote also by $\alpha_0$ the contact form on $S^*\Funk^\eta_L(K)$, and by $\overline\alpha_0$ that on $S^*\Funk^\eta{K^\vee}(L^\vee)$. 
Define $\alpha_\eta =(\Phi_\eta^{-1})^*\alpha_1\in\Omega^1(\Image(\Phi_\eta))$ and $\overline\alpha_\theta =(\overline\Phi_\theta^{-1})^*\alpha_1\in\Omega^1(\Image(\Phi_\eta))$, which coincide with the restriction to $\Image(\Phi_\eta)$ of the corresponding $1$-forms on $\partial K\times\partial L^\vee$ in the proof of Theorem \ref{thm:dual_lengths}, in particular by that proof $\alpha_\eta-\overline\alpha_\theta$ is an exact form.

We now have

 \begin{align*}\vol^{HT}(\Funk_L(\partial K))&=\int_{S^*\partial K}(d\alpha_1)^{n-2}\wedge\alpha_1=\int_{\Image(\Phi_\eta)}d\alpha_\eta^{n-2}\wedge\alpha_\eta\\&=\int_{\Image(\overline\Phi_\theta)}d\overline\alpha_\theta^{n-2}\wedge\overline\alpha_\theta=\vol^{HT}(\Funk_{K^\vee}(\partial L^\vee)).\end{align*}
 Similarly, the length spectrum of $\Funk^\eta_L(\partial K)$ coincides with the action spectrum of $\alpha_\eta$ on $\Image(\Phi_\eta)$, which is also the action spectrum of $\overline \alpha_\theta$.
%
%Then $\Phi_\eta^*(d\alpha_\eta-\omega)=0\Rightarrow d\alpha_\eta|_{\Image(\Phi_\eta)}=\omega|_{\Image(\Phi_\eta)}$, and so $\alpha_\eta|_{\Image(\Phi_\eta)}$ is independent of $\eta$ up to the addition of an exact 1-form.
%
%Since $$\vol^{HT}(\Funk^\eta_L(\partial K))=\int_{S^*\partial K}(d\alpha_1)^{n-2}\wedge\alpha_1=\int_{\Image(\Phi_\eta)}\omega^{n-2}\wedge\alpha_\eta$$ is unchanged when $\alpha_\eta$ is perturbed by an exact form, while $\Image(\Phi_\eta)$ is independent of $\eta$, the projective invariance of the volume follows. Similarly, the length spectrum coincides with the action spectrum of $\alpha_\eta$ on $\Image(\Phi_\eta)$, which is again independent of perturbations by exact forms.
%
%Turning to the duality statement, assume $K$ is strictly convex and define $\overline{\Phi_\theta}=\overline{H_\theta}^{-1}\circ\tau_\theta:S^*\partial L^\vee\to \partial K\times\partial L^\vee$ in full analogy with $\Phi_\eta$. Observe that the images of the two maps coincide, since the degeneracy of the pairing  $T_x\partial K\otimes T_\xi\partial L^\vee\to \R$ is symmetric in $(x,\xi)$. As $\overline{\Phi_\theta}^*\omega=\Phi_\eta^*\omega=d\alpha_1$, and in light of the proof of projective invariance, the statement follows.
\endproof
}
We remark that the Holmes-Thompson volumes of $\Funk_L(\partial K)$ and $\RFunk_L(\partial K)$ trivially coincide.

\textit{Proof of Corollary \ref{cor:hilbert_plane}.}
This follows at once from Theorem \ref{thm:girth}, since the Hilbert length of a curve coincides with its 1-dimensional Funk Holmes-Thompson volume.
\qed.

\section{Mahler volumes of every Funk scale}\label{sec:every_scale}
Define the \emph{Funk-Mahler volume of radius $0<r<\infty$ (around $q\in \inter(K)$)} by $$M_r(K,q)=\omega_n r^{-n}\vol_K^F(B_r^F(q, K)), \quad  M_r(K)=\inf_{q\in \inter(K)}M_r(K,q).$$

In Euclidean terms, we have $$M_r(K,q)=r^{-n}\int_{z\in q+(1-e^{-r})(K-q)}|K^z|dz,$$  where $K^z$ denotes the dual body with respect to $z$. Evidently, $M_r(K,q)$ is invariant under the full linear group with $q$ at the origin. 

\begin{Lemma}\label{lem:mahler_convex}
	For any convex $K$, $M_r(K, q)$ is a strictly convex function of $q$. In particular, if $K=-K$ then $M_r(K)=M_r(K, 0)$.
\end{Lemma}
\proof
The function $z\mapsto |K^z|$ is well-known to be strictly convex in $z\in \inter(K)$. This follows at once by writing 
$$|K^z|=\omega_n\int_{S^{n-1}}\max_{x\in K}\langle \theta, x-z\rangle ^{-n}d\sigma(\theta),$$
where $\sigma$ is the rotation-invariant probability measure on $S^{n-1}$, and noticing that for any fixed $x,\theta$, the function $z\mapsto (\langle x,\theta\rangle -\langle z, \theta\rangle)^{-n}$ is strictly convex.

It follows that

$$M_r(K,q)=r^{-n}\int_{(1-e^{-r})K}|K^{e^{-r}q+w}|dw$$ is convex as well.

If $K=-K$ then $M_r(K, q)=M_r(K, -q)$, and so $q=0$ must be a minimum.

\endproof

The Funk-Mahler volume enjoys invariance under duality, as follows.

\begin{Proposition}
	It holds for all $K\subset \R^n$ that $M_r(K, q)=M_r(K^q, q)$. In particular, if $K=-K$ then $M_r(K)=M_r(K^o)$.
\end{Proposition}
\proof
We may assume $q=0$. Put $\rho=1-e^{-r}$, so that $B_r^F=\rho K$ is the ball of radius $r$ in the Funk metric, centered at the origin. Then  $$M_r(K,0)=\frac{\omega_n}{r^n}\vol_K^F(\rho K)=\frac{\omega_n}{r^n}\vol_{\rho^{-1} K^o}^F(K^o)=M_r(\rho^{-1}K^o,0)=M_r(K^o, 0),$$ where the second equality is due to Corollary \ref{cor:dual_volumes}, the penultimate equality holds since $K^o$ is the ball of radius $r$ around the origin in the Funk geometry of $\rho^{-1}K^o$, and the last equality is due to $\GL$-invariance.

The second statement now follows by Lemma \ref{lem:mahler_convex}.
\endproof
\begin{Remark}
In projective-invariant terms, this duality assumes the following form. For $K\subset\PP$, $q\in \inter(K)$ and $\eta\in\inter(K^\vee)$, let $B^\eta_r(K, q)$ be the outward ball of radius $r$ centered at $q$ in $\Funk_\eta(K)$. Then $B^\eta_r(K, q)\subset \Funk_\eta(K)$ and $B_r^q(K^\vee, \eta)\subset \Funk_q(K^\vee)$ have equal volume.
\end{Remark}

As $r\to0$, $M_r(K,q)\to M_0(K,q):=|K||K^q|$. Thus the Mahler volume $M_0(K):=\lim_{r\to 0} M_r(K)$ is recovered at infinitesimal scale.

On the opposite end of the scale, the asymptotics of $M_r(K,q)$ is governed by the regularity of the boundary of $K$. When $K$ is smooth and strictly convex, we recover the centro-affine surface area as $r\to\infty$, as established in \cite{berck_bernig_vernicos}. Recall that the \emph{centro-affine surface area relative to $q\in \inter(K)$} is given by 

\begin{equation}\label{eq:centro_affine}\Omega_n(K, q)=\int_{\partial K}\frac{k_p^{\frac 12}}{\langle p-q, \nu_p\rangle^{\frac{n+1}{2}}}d\mu_K(p;q)=\int_{\partial K}\frac{k_p^{\frac 12}}{\langle p-q, \nu_p\rangle^{\frac{n-1}{2}}}d\mathcal H^{n-1}(p),\end{equation}
where $k_p$ is the Gaussian curvature of $\partial K$ at $p$, $\nu_p$ the unit outward normal, and $\mu_K(p;q)$ the cone measure with the origin at $q$.

As \cite{berck_bernig_vernicos} focuses on the Hilbert metric with the Busemann volume, we provide for the sake of completenss the details of the computation in our case.

\begin{Lemma}\label{lem:volume_ellipsoid}
	Let $A$ be a positive-definite $n\times n$ matrix, and $Q_A=\{x\in\R^n:\langle Ax, x\rangle \leq 1\}$. Then for $p\in\partial Q_A$ and $0<\rho <1$, we have 
	$$|Q_A^{\rho p}|=\omega_n\frac{(\det A)^{1/2}}{(1-\rho^2)^{\frac{n+1}{2}}}=\omega_n\frac{1}{(1-\rho^2)^{\frac{n+1}{2}}}\frac{k_p^{1/2}}{\langle p,\nu_p\rangle^{\frac{n+1}{2}}}.$$ 
\end{Lemma}
\proof
For the first equality, we compute that $Q_A^{\rho p}=\{ \xi: \langle B(\xi-c),\xi-c\rangle\leq 1\}$, where, denoting $B_0=A^{-1}-\rho^2pp^T$, we have 
\begin{align*}c&=\rho B_0^{-1}p, \qquad B=\frac{B_0}{1+\rho^2\langle B_0^{-1}p,p\rangle}.
\end{align*}
Now by the Sherman-Morrison formula, $$B_0^{-1}=A+\frac{\rho^2}{1-\rho^2}App^TA\Rightarrow  \langle B_0^{-1}p,p\rangle =\frac{1}{1-\rho^2},$$ while by the matrix determinant lemma we h	ave
$$\det B_0 =(1-\rho^2)\det A^{-1}.$$

It follows that
$$|Q_A^{\rho p}|=\omega_n\sqrt{\det B}^{-1}=\omega_n\sqrt{\frac{\det A}{1-\rho^2}}(1-\rho^2)^{-\frac{n}{2}}=\omega_n\frac{\sqrt{\det A}}{(1-\rho^2)^{\frac{n+1}{2}}}$$

For the second equality, recall that $\nu_x=\frac{Ax}{|Ax|}$. It follows that for $v\in T_pQ_A$, $$D_p\nu_x (v)=\frac{Av}{|Ap|}-\langle \frac{Av}{|Ap|}, \nu_p\rangle \nu_p=\frac{1}{|Ap|}\pi_{T_pQ_A}\circ A.$$

Take an orthonormal basis $v_1,\dots, v_{n-1}$ of $T_pQ_A$, and evaluate the volume of $A(v_1\wedge\dots\wedge v_{n-1}\wedge p)$ in two ways. We get $$|A(v_1\wedge\dots\wedge v_{n-1}\wedge p)|=\det A\cdot\langle p, \nu_p\rangle=\det(\pi_{T_pQ_A}\circ A)\cdot|Ap|.$$
Noting also that $\langle p, \nu_p\rangle=\frac{1}{|Ap|}$ since $\langle Ap, p\rangle =1$, we find
$\det A=\frac{k_p}{\langle p, \nu_p\rangle^{n+1}}.$
\endproof

\begin{Lemma}\label{lem:dual_volume}
	For $K\in\mathcal K_0(\R^n)$ smooth and strictly convex, and $p\in \partial K$, we have $$|K^{\rho p}|=  (1-\rho)^{-\frac{n+1}{2}}\left(2^{-\frac{n+1}{2}} \omega_n \frac{k_p^{1/2}}{\langle p,\nu_p\rangle^{\frac{n+1}{2}}}+err(p)\right),$$
	where $err(p)\to 0$ as $\rho\to 1$, uniformly in $p\in\partial K$.
\end{Lemma}
\proof
Throughout the proof, $c,c'$ will be positive constants that are bounded away from $0$ and $\infty$, uniformly in $p\in\partial K$, and may change between occurrences.

As both $|K^{\rho p}|$ and the main term on the right hand side are $\SL(n)$-invariant, we may assume without loss of generality that $p=-e_n$ and $T_p\partial K=\Span(e_1,\dots,e_{n-1})$. We may further assume that $y_{n}=-1+\frac {a_p} 2|y'|^2$ is the quadratic approximation of $\partial K$ near $p$, where we write $y'=(y_1,\dots,y_{n-1})$. Let $Q$ be the ellipsoid $a_p|y'|^2 +y_n^2=1$, osculating $\partial K$ at $p$. Write $1-\rho =\epsilon^2$, and note that $c\leq a_p\leq c'$.

Let $\nu_x$ be the Euclidean normal at $x\in\partial K$. One has 
\[\|\nu_x\|_{K^{\rho p}}=h_K(\nu_x)-\langle \nu_x,\rho p\rangle=\langle x-\rho p, \nu_x\rangle. \]
We may therefore parametrize the points of $K^{\rho p}$ by $$\zeta_K:[0,1]\times\partial K\to K^{\rho p}, \quad \zeta_K(\alpha, x)=\alpha\frac{\nu_x}{\langle x-\rho p, \nu_x\rangle}.$$
 Then by Cauchy-Schwartz,
$$\langle x-\rho p, \nu_x\rangle\geq \langle x-p, \nu_x\rangle -(1-\rho)|\langle p, \nu_x\rangle| \geq  \langle x-p, \nu_x\rangle- \epsilon^2.$$

Let $H$ be the hyperplane parallel to $T_p\partial K$ at distance $\epsilon$ from it, and $H^+, H^-$ the half-spaces it defines, with $p\in H^+$. We claim that $|\zeta_K([0,1]\times \partial (K\cap H^-))|\leq c\epsilon^{-n}	$.

Consider $x\in\partial K\cap H^-$, so $|x-p|\geq \epsilon$. If $\nu_x$ is proportional to $x-p$ then $\langle x-p, \nu_x\rangle |x-p|\geq \epsilon$, and $|\zeta_K(\alpha, x)| \leq \frac2 {\epsilon}$. Otherwise, consider the affine plane $P=p+\Span(x-p, \nu_x)$, and consider the planar convex body $K_P=K\cap P$. Using Euclidean cooridnates $(x_1,x_2)$ in which $T_pK_P=\R e_1$, $Q_P=Q\cap P$ is osculated at $p$ by the parabola $x_2=-1+\frac12 a_p x_1^2$.

\begin{figure}[h]
	\centering
	\includegraphics[width=0.4\linewidth]{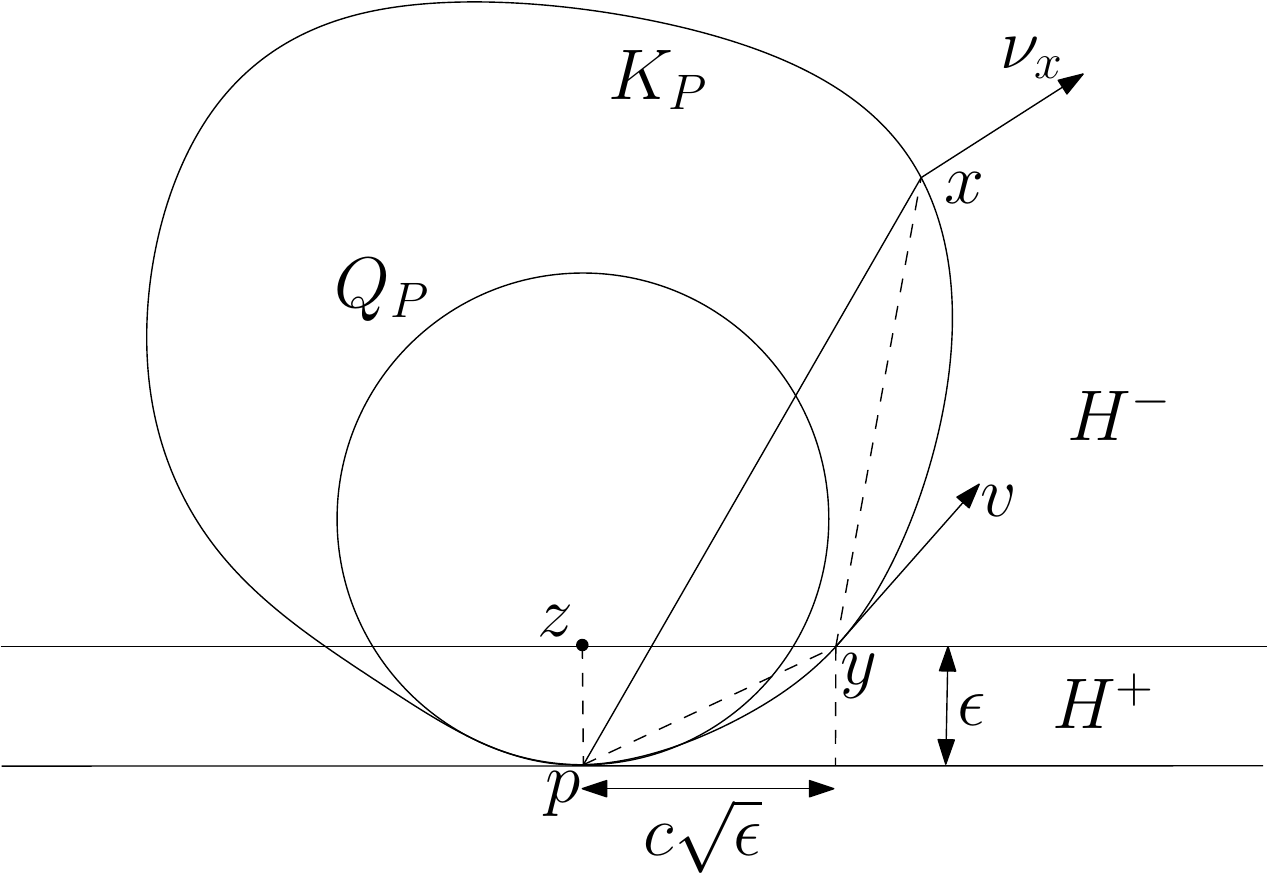}
	\caption{}
	\label{fig:reflectionlaw}
\end{figure}

We orient $\partial K$ in such a way that the tangent at $x$ makes an acute angle with $x-p$. Let $y\in H\cap \partial K_P$ belong to the oriented arc $[p,x]\subset\partial K_P$.  
It holds by the convexity of $K_P$ that $$\langle x-p, \nu_x\rangle\geq |x-p|\sin \measuredangle pxy=|p-y|\sin\measuredangle pyx\geq c\sqrt{\epsilon}\sin\measuredangle pyx.$$
Denoting by $v$ the tangent ray to $\partial K_P$ at $y$, we have $$c\sqrt \epsilon=\measuredangle pyz\leq \measuredangle pyx\leq \pi-\measuredangle(y-p, v)= \pi-c\sqrt \epsilon.$$ Therefore, $$\langle x-p, \nu_x\rangle\geq c\sqrt\epsilon  \sin\measuredangle pyx\geq c\epsilon.$$

Thus again we find that $|\zeta_K(\alpha, x)|\leq \frac c\epsilon$. That is, $\zeta_K([0,1]\times \partial K\cap H^-)\subset \frac{c}{\epsilon}B^n$, and so \begin{equation}\label{eq:volume_bound}|\zeta_K([0,1]\times(\partial K\cap H^-))|\leq c\epsilon^{-n}.\end{equation}

Define $z_K(y'), z_Q(y'):\R^{n-1}\to \R$ for $|y'|<c\sqrt{\epsilon}$ by having $x_K(y'):=(y', z_K(y'))\in \partial K, x_Q(y'):=(y', z_Q(y'))\in \partial Q$, and let $\nu_K(y'), \nu_Q(y')$ be the correspondning unit normals. Then $z_K(y')=z_Q(y')+O(|y'|^3)=-1+\frac{a_p}{2}|y'|^2+O(|y'|^3)$.

\begin{figure}[h]
	\centering
	\includegraphics[width=0.4\linewidth]{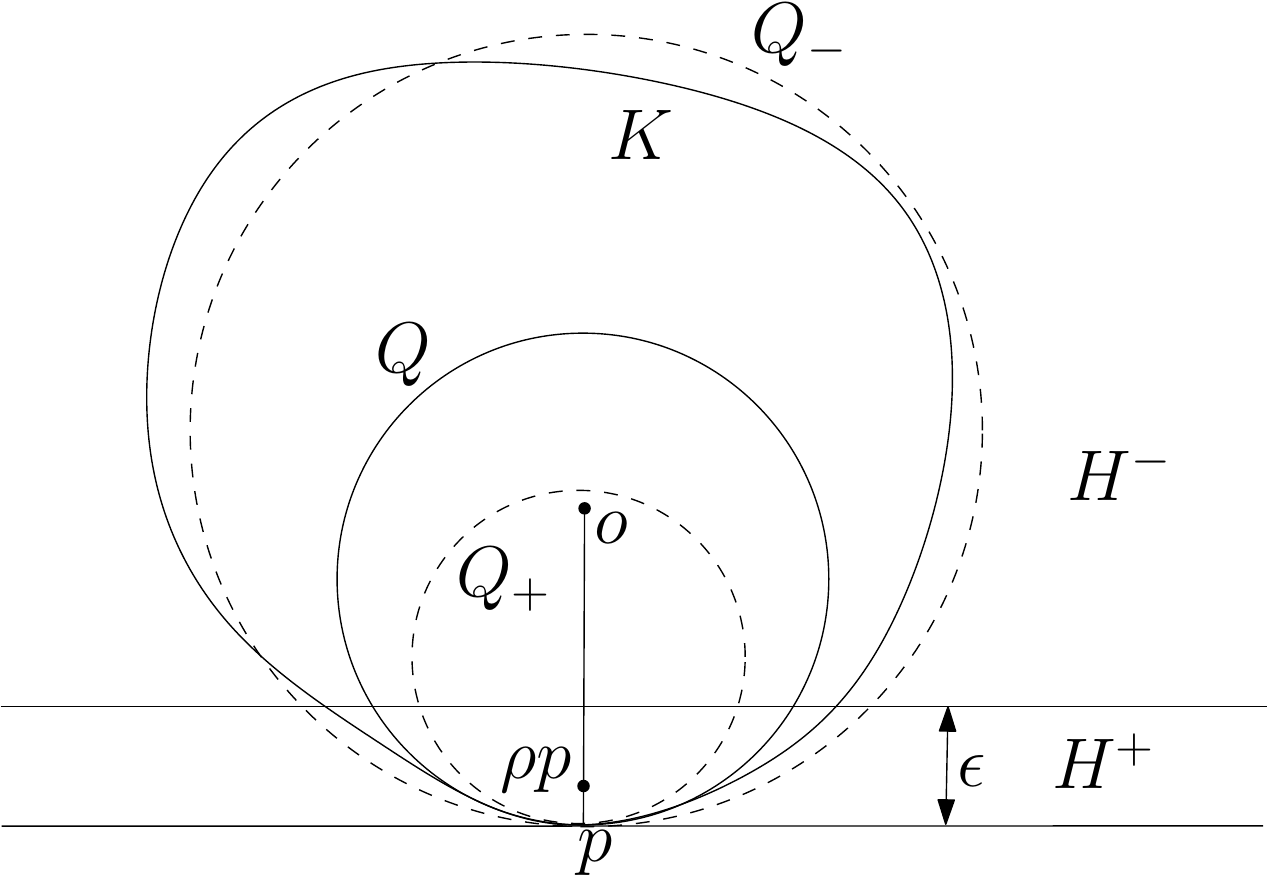}
	\caption{}
	\label{fig:reflectionlaw}
\end{figure}

Define the ellipsoids $Q_{\pm} = \{y: \frac12 (a_p\pm 2\epsilon^{1/3})|y'|^2+y_n^2=1\}.$
It holds that $z_{Q_\pm}(y')=-1+\frac{1}{2}(a_p\pm 2\epsilon^{1/3})|y'|^2+O(|y'|^3),$
and so for $\epsilon\ll1$, $$Q_+\cap H^+\subset K\cap H^+\subset Q_-\cap H^+.$$
Moreover, as $|\nu_K(y')-\nu_Q(y')|=O(|y'|)$ and $|y'|\leq c \epsilon^{1/2}\ll\epsilon^{1/3}$, it is evident that $$\zeta_{Q_-}([0,1]\times Q_-\cap H^+)\subset  \zeta_{K}([0,1]\times K\cap H^+)\subset \zeta_{Q_+}([0,1]\times Q_+\cap H^+).$$
By estimate \eqref{eq:volume_bound} applied to $Q_\pm$ and combined with Lemma \ref{lem:volume_ellipsoid}, we have
$$ \zeta_{Q_\pm}([0,1]\times Q_\pm\cap H^+)=\epsilon^{-(n+1)}\left( 2^{-\frac{n+1}{2}} \omega_n \frac{k_p^{1/2}}{\langle p,\nu_p\rangle^{\frac{n+1}{2}}}+err_\pm(p)\right), $$
and applying \eqref{eq:volume_bound} for $K$, the claim follows.

\endproof

\begin{Proposition}\label{prop:centro_affine}
	For $q\in \inter(K)$, $$\lim_{R\to \infty}\frac{\vol^F_K(B_R^F(q, K))}{e^{\frac{n-1}{2}R}}=\frac{1}{n-1}2^{-\frac{n-1}{2}} \Omega_n(K,q).$$
\end{Proposition}
\proof
We may put $q=0$. Parametrize the outward Funk ball of radius $R$ around $0$ by $p\in\partial K$, $0\leq r\leq R$ using $\phi(p,r)=(1-e^{-r})p$. Then

$$ \vol^F_K(B_R^F(0,K))=\int_0^R\int_{\partial K} |K^{\phi(p,r)}|\mathrm{Jac}(\phi(p,r))d\mathcal H^{n-1}(p)dr.$$  
Now $\mathrm{Jac}(\phi(p,r))=(1-e^{-r})^{n-1}e^{-r}\langle p, \nu_p\rangle$, while by Lemma \ref{lem:dual_volume} 
$$ |K^{\phi(p,r)}|= 2^{-\frac{n+1}{2}} \omega_n e^{\frac{n+1}{2}r} \frac{k_p^{1/2}}{\langle p,\nu_p\rangle^{\frac{n+1}{2}}}+o(e^{\frac{n+1}{2}r}),$$
where the error term is uniform in $p$. Applying l'Hospital's rule we find 
\begin{align*}\lim_{r\to \infty}\frac{\vol_K^F(B_r^F(0,K))}{e^{\frac{n-1}{2}r}}&=\lim_{r\to\infty}\frac{2}{n-1} \frac {(1-e^{-r})^{n-1}e^{-r}\int_{\partial K} |K^{\phi(p,r)}|\langle p, \nu_p\rangle  d\mathcal H^{n-1}(p)}{e^{\frac{n-1}{2}r}}
\\&=\frac{1}{n-1}2^{-\frac{n-1}{2}} \int_{\partial K}\frac{k_p^{1/2}}{\langle p,\nu_p\rangle^{\frac{n-1}{2}}}d\mathcal H^{n-1}(p).\end{align*}
\endproof

\section{Maximizing the Funk-Mahler volume}\label{sec:inequalities}

Let us introduce for convenience $\rho=1-e^{-r}$, and the modified Funk-Mahler volumes
$$\widetilde M_\rho(K,q):=r^nM_{r}(K,q)=\int_{z\in q+\rho (K-q)}|K^z|dz,\quad \widetilde M_\rho(K)=\inf_{q\in \inter(K)}\widetilde M_\rho(K, q).$$

\begin{Lemma}
	Let $B\subset\R^n$ be an ellipsoid. Then $$\widetilde M_\rho(B)=n\omega_n^2\int_0^\rho\frac{t^{n-1}dt}{(1-t^2)^{\frac{n+1}{2}}}.$$
\end{Lemma}
\proof By Lemmas \ref{lem:mahler_convex} and \ref{lem:volume_ellipsoid},
\[\widetilde M_\rho(B)=\widetilde M_\rho(B^n,0)=\int_{\rho B}|B^z|dz=\omega_n\int_{\rho B}\frac{dz}{(1-|z|^2)^{\frac{n+1}{2}}}=n\omega_n^2  \int_0^\rho\frac{t^{n-1}dt}{(1-t^2)^{\frac{n+1}{2}}}.\]
\endproof

In light of the discussion in section \ref{sec:every_scale}, and inspired by the Blaschke-Santal\'o and the centro-affine isoperimetric inequalities, we propose a unified inequality containing both of the above as limiting cases.
\begin{Conjecture}\label{conj:bodies}
	Fix any $0\leq \rho<1$, and let $K\in\mathcal K(\R^n)$ have $0$ as its centroid.
	\begin{enumerate}
		\item $\widetilde M_\rho(K,0)\leq \widetilde M_\rho(B^n)$.
		\item If $\widetilde M_\rho(K,0)= \widetilde M_\rho(B^n)$ then $K$ is an ellipsoid.
	\end{enumerate}
\end{Conjecture}
This is not to say that $\widetilde M_\rho(K)=\widetilde M_\rho(K,0)$: as $\rho\to 0$, $\widetilde M_\rho(K,q)$ attains its minimum close to the Santal\'o point of $K$, which generally differs from the centroid.

We will prove the conjecture in the restricted class of unconditional convex bodies. To this end, we must expand our discussion to include log-concave functions. We write $f_q(\bullet)=f(q+\bullet)$.

\begin{Definition}
	For $f:\R^n\to \R_+$ a Borel function, put $\phi(x)=-\log f(x):\R^n\to \R\cup\{+\infty\}$.  Define for $\rho\in\mathbb C$ the Funk-Mahler volumes around $0$
	\begin{align*}\widetilde M_\rho(f,0)&:=\frac{1}{n!}\displaystyle\int_{z\in\R^n} e^{-\phi(\frac{z}{\rho})}\displaystyle\int_{\xi\in \R^n} e^{-\mathcal L\phi_z(\xi)} d\xi dz\\&=\frac{1}{n!}\rho^n\displaystyle\int_{\R^n\times \R^n} e^{-\phi(x)-\mathcal L \phi (\xi)+\rho\langle x,\xi\rangle}dxd\xi,\end{align*}
	whenever the integral converges. Define also
	$$\widetilde M_\rho(f, q):=\widetilde M_\rho(f_q,0), \quad \widetilde M_\rho(f):=\inf_{q\in \R^n} \widetilde M_\rho(f, q).$$
\end{Definition}
As $\rho\to 0$, $$\widetilde M_\rho(f,0)\sim \frac 1{n!} \rho^n \int_{\R^n\times \R^n} e^{-\phi(z)-\mathcal L \phi(\xi)}dxd\xi,$$ which is the functional Mahler product.

 \begin{Lemma}\label{lem:functional_invariance}
 It holds that $\widetilde M_\rho(f,0)=\widetilde M_\rho(\lambda f\circ A,0)$ for any $\lambda> 0$ and $A\in\GL_n(\R)$. Furthermore, $\widetilde M_\rho(f,q)$ is a convex function of $q$, and strictly convex where finite.
 \end{Lemma}
 \proof
 The first statement is immediate from definition by the $\GL$-equivariance of the Legendre transform, and since $\mathcal L(\phi+c)=\mathcal L\phi -c$.

 For the second, observe that
 $$ \widetilde M_\rho(f,q)=\frac{1}{n!}\int_{\R^n\times\R^n} e^{-\phi(\frac{w}{\rho})+\langle \xi,w\rangle-\mathcal L\phi(\xi)}e^{(1-\rho)\langle \xi,q\rangle}dwd\xi,$$
and note that $t\mapsto e^{(1-\rho)t}$ is strictly convex.
 \endproof

	The functional Funk-Mahler volumes generalize those previously introduced for convex bodies, as follows.
	
	\begin{Lemma}\label{lem:body_vs_function}
		$\widetilde M_\rho(\mathbbm 1_K, 0)=\widetilde M_\rho(K ,0)$.
	\end{Lemma} 
	\proof
	Recall that $\mathcal L(-\log \mathbbm 1_K)(\xi)=h_K(\xi)=\|\xi\|_{K^o}$. Thus 
	$$\widetilde M_\rho(\mathbbm 1_K, 0)=\frac{1}{n!}\int_{\rho K\times\R^n}e^{-\|\xi\|_{K^z}}dzd\xi=\frac{1}{n!}\int_{\rho K}dz \int_0^\infty r^{n-1}e^{-r}dr\int_{\partial K^z}d\mu_{K^z}, $$
	where $d\xi=r^{n-1}drd\mu_{K^z}$. Thus 
	$$\widetilde M_\rho(\mathbbm 1_K, 0)= \frac{1}{n!}(n-1)!n\int_{\rho K} |K^z|dz=\widetilde M_\rho(K, 0). \qedhere$$

   Using this representation, let us exhibit a curious property of the Funk-Mahler volume of a convex body, which at present we only formally establish for zonoids. 
 
 \begin{Proposition}\label{prop:to_the_infinity_and_beyond} Assume $0\in \inter(K)$. There is then a continuous extension of $\widetilde M_\rho(K, 0)$ to $\rho\in[0,1)\times i\R$, which is analytic in $(0,1)\times i\R$. 
  \\Furthermore if $K$ is a zonoid, then
 	$$\lim_{t\to \pm\infty} \widetilde M_{it}(K,0)=i^n\frac{(2\pi)^n}{n!}.$$
 \end{Proposition}
 \proof
 
 For the first statement, put $\rho=s+it$, write 
  $$\widetilde M_{\rho}(K,0)= \widetilde M_{\rho}(\mathbbm 1_K,0)=\frac{(s+it)^n}{n!}\int_{K\times\R^n} e^{-h_K(\xi)+s\langle x,\xi\rangle +it \langle x,\xi\rangle}dxd\xi, $$
   and note that 
  
 $$\frac{1}{n!}\int_{K\times\R^n} e^{-h_K(\xi)+s\langle x,\xi\rangle }dxd\xi =\frac{1}{s^n}\widetilde M_s(K,0)<\infty$$
for $0\leq s<1$.

For the second statement, let $\mathcal F(f)(y)=\int_{\R^n}f(x)e^{-i\langle x, y\rangle }dx$ be the Fourier transform. Write 
\begin{align*}\widetilde M_{it}(K,0)&=\frac{i^n}{n!} t^n\int_{\R^n\times\R^n}\mathbbm 1_K(x) e^{-h_K(\xi)}e^{it\langle x, \xi\rangle}dxd\xi\\&=\frac{i^n}{n!}\int_{\R^n} \mathbbm 1_K(\frac{x}{t})\mathcal F(e^{-h_K})(-x)dx\\&=\frac{i^n}{n!}\int_{\R^n}e^{-h_K(\xi)}\int_{tK}e^{i\langle x,\xi\rangle} dxd\xi ,\end{align*}
and note that $\mathcal F(e^{-h_K})\in L^1(\R^n)$ when $K$ is a zonoid. Thus as $t\to\pm \infty$, 
$$ \widetilde M_{it}(K,0)\to \frac{i^n}{n!}\int_{\R^n}\mathcal F(e^{-h_K})=\frac{i^n}{n!}(2\pi)^ne^{-h_K(0)}=i^n\frac{(2\pi)^n}{n!}.$$
 \endproof
 
Motivated by Conjecture \ref{conj:bodies}, and inspired by the results of \cite{lehec}, we boldly put forward the following
\begin{Conjecture}\label{conj:functional_general}
	Consider two Borel functions $\phi,\psi:\R^n\to \R\cup\{+\infty\}$ such that for all $x,\xi\in\R^n$ one has $\phi(x)+\psi(\xi)\geq \langle x, \xi\rangle $. Assume also $\int xe^{-\phi(x)}dx=0$. Then for all $0\leq \rho<1$, 
		\begin{equation}\label{eq:general_inequality} \int_{\R^n\times\R^n}e^{-\phi(x)-\psi(\xi)+\rho\langle x,\xi\rangle}dxd\xi\leq (2\pi)^n(1-\rho^2)^{-\frac n 2}.\end{equation}
	Equality is attained uniquely, up to equality a.e., by the pairs of functions
	\begin{equation}\label{eq:pairs}\phi(x)=\frac 12 \langle Ax, x\rangle +c, \quad \psi(\xi)=\mathcal L\phi(\xi)=\frac 12 \langle A^{-1}x, x\rangle -c,\end{equation}
	for some $A>0$ and $c\in\R$.
\end{Conjecture}
We remark that here and in the following, one may restrict attention to dual pairs of convex functions $\phi, \psi=\mathcal L\phi$, without loss of generality.
In all that follows, we take \emph{a.e.-uniqueness} to mean uniqueness up to equality a.e.

Let us check the trivial statement in the conjecture. 
\begin{Lemma}
	Pairs $\phi,\psi$ as in \eqref{eq:pairs} give equality in \eqref{eq:general_inequality}.
\end{Lemma}
\proof 
Taking $q(x)=\frac12 |x|^2$, $$\mathcal Lq_z(\xi)=\frac{|\xi|^2}{2}-\langle \xi,z\rangle,$$ 
and so
\begin{align*}
&\int_{\R^n\times\R^n} e^{-\frac{|z|^2}{2\rho^2}-\frac{|\xi|^2}{2}+\langle \xi, z\rangle}d\xi dz=\rho^n\int e^{-\frac12 (|w|^2+|\xi|^2-2\rho\langle w,\xi\rangle )}d\xi dw
	\\&=\rho^n\int e^{-\frac12 (|w-\rho\xi|^2+(1-\rho^2)|\xi|^2)}d\xi dw =\int e^{-\frac12 |w|^2}dw\int e^{-\frac12 (1-\rho^2)|\xi|^2} d\xi \\&=(2\pi)^n\rho^n(1-\rho^2)^{-\frac{n}{2}}.
\end{align*}
By Lemma \ref{lem:functional_invariance}, we get equality for all pairs as claimed.

\endproof
An immediate corollary is a functional version of Conjecture \ref{conj:bodies}, as follows.
\begin{Conjecture}\label{conj:functional_santalo}
		Let $f:\R^n\to \R_+$ be an integrable Borel function, and assume also that $\int x f(x)dx=0$. Then for all $0\leq \rho<1$,
		\[ \widetilde M_\rho(f,0)\leq \frac{1}{n!}(2\pi)^n\rho^n(1-\rho^2)^{-\frac n 2}. \]
		Equality is attained uniquely by multiples of gaussians: $f(x)=\lambda e^{-\frac 12 \langle Ax, x\rangle}$.
\end{Conjecture}
\begin{Remark}
Note that the implication Conjecture \ref{conj:functional_general} $\Rightarrow$ Conjecture \ref{conj:functional_santalo} is valid separately for each $n$, and also separately for any family of functions with certain symmetry, such as even or unconditional functions.
\end{Remark}

We now prove Conjecture \ref{conj:functional_general} under the extra assumption of unconditionality. Recall that $\phi$ is unconditional if $\phi(x_1,\dots, x_n)=\phi(|x_1|,\dots, |x_n|)$.

\begin{Theorem}\label{thm:functional_unconditional}
	Consider two unconditional Borel functions $\phi,\psi:\R^n\to \R\cup\{+\infty\}$ such that for all $x,\xi\in\R^n$ one has $\phi(x)+\psi(\xi)\geq \langle x, \xi\rangle $. Then for all $0\leq \rho<1$, $$ \int_{\R^n\times\R^n}e^{-\phi(x)-\psi(\xi)+\rho\langle x,\xi\rangle}dxd\xi\leq (2\pi)^n(1-\rho^2)^{-\frac n 2}.$$
	Equality is attained a.e.-uniquely by Legendre-dual pairs $$\phi(x)=\frac 12 \langle Ax, x\rangle +c, \quad \psi(\xi)=\mathcal L\phi(\xi)=\frac 12 \langle A^{-1}x, x\rangle -c.$$
 \end{Theorem}

\proof

Since $\phi(x)=\phi(-x)$, the claimed inequality is equivalent to 
\begin{equation}\label{eq:ineq1} I:=\int_{\R^n\times\R^n} e^{-\phi(x)-\psi(\xi)}\cosh(\rho \langle x, \xi\rangle)d\xi dx\leq (2\pi)^n(1-\rho^2)^{-n/2}.\end{equation}

The left hand side can be rewritten as 
\begin{align*}J&=\sum_{j=0}^\infty \frac{\rho^{2j}}{(2j)!}\int_{\R^n\times\R^n}e^{-\phi(x)-\psi(\xi)}\langle x, \xi\rangle^{2j}dxd\xi\\&=\sum_{j=0}^\infty\frac{\rho^{2j}}{(2j)!}\sum_{|I|=2j}{2j\choose I}\int_{\R^n}x^{I}e^{-\phi(x)}dx\int_{\R^n}\xi^{I}e^{-\psi(\xi)}d\xi,\nonumber\end{align*}
where the last sum is over all $n$-tuples $I=(i_1,\dots, i_n)$ of non-negative integers such that $|I|:=\sum_{\nu=1}^n i_\nu =2j$, and we use the notation $x^I:=\prod_{\nu=1}^n x_\nu^{i_\nu}$, and $${2j\choose I}:=\frac{(2j)!}{i_1!\cdots i_n!}$$ for the multinomial coefficient. 

Observe that by unconditionality, if $I$ contains an odd index then the corresponding integrals vanish. Therefore, 
\begin{align}\label{eq:eq1} J&= \sum_{j=0}^\infty\frac{\rho^{2j}}{(2j)!}\sum_{|I|=j}{2j\choose 2I}\int_{\R^n}x^{2I}e^{-\phi(x)}dx\int_{\R^n}\xi^{2I}e^{-\psi(\xi)}d\xi\\&=
	4^{n}\sum_{j=0}^\infty\frac{\rho^{2j}}{(2j)!}\sum_{|I|=j}{2j\choose 2I}\int_{\R^n_+}x^{2I}e^{-\phi(x)}dx\int_{\R^n_+}\xi^{2I}e^{-\psi(\xi)}d\xi. \nonumber\end{align}

Making the change of coordinates $x=e^u:=(e^{u_i})_{i=1}^n$, $\xi=e^v=(e^{v_i})$, we may write 
\begin{equation}\label{eq:eq2}\int_{\R_+^n}e^{-\phi(x)}x^{2I}dx\int_{\R_+^n}e^{-\psi(\xi)}\xi^{2I}d\xi=\int_{\R^n} f(u)du\int_{\R^n} g(v)dv,\end{equation}
where $f(u)=e^{-\phi(e^u)}e^{\langle 2I+1, u\rangle}$, $g(v)=e^{-\psi(e^v)}e^{\langle 2I+1, v\rangle}$, and $I+c:=(i_\nu+c)_{\nu=1}^n$.

Denote $$h(w)=e^{-\frac12 |e^w|^2}e^{\langle 2I+1, w\rangle},\quad w\in \R^n,$$ where $|e^w|=\sqrt{\sum_{i=1}^n e^{2w_i}}$. Since $\phi(x)+\psi(\xi)\geq  \langle x, \xi\rangle$, we immediately find that
$$f(u)g(v)\leq h\left(\frac{u+v}{2}\right)^2,$$
and by the Prekopa-Leindler inequality we have
\begin{align*} \int_{\R^n} f\int_{\R^n} g&\leq \left(\int_{\R^n} h\right)^2=\left(\int_0^\infty e^{-\frac{|p|^2}{2}}p^{2I}dp\right)^2=\prod_{\nu=1}^n\left(\int_0^\infty e^{-\frac{t^2}{2}} t^{2i_\nu}dt\right)^2\\&=\left( \frac{\sqrt{2\pi}}{2}\right)^{2n} \prod_{\nu=1}^n ((2i_\nu-1)!!)^2=\frac{(2\pi)^n}{4^n} \prod_{\nu=1}^n ((2i_\nu-1)!!)^2.\end{align*}
Combined with \eqref{eq:ineq1}, \eqref{eq:eq1}, \eqref{eq:eq2}, it remains to verify that

$$ \sum_{j=0}^\infty \frac{\rho^{2j}}{(2j)!}\sum_{|I|=j}{2j\choose 2I}\prod_{\nu=1}^n ((2i_\nu-1)!!)^2 =(1-\rho^2)^{-n/2},$$
or equivalently 

$$ \sum_{j=0}^\infty \rho^{2j}\sum_{|I|=j}{2j\choose 2I}\prod_{\nu=1}^n \frac{(2i_\nu-1)!!}{(2i_\nu)!!} =(1-\rho^2)^{-n/2}.$$
By Newton's binomial formula, 
\[ (1-\rho^2)^{-\frac12}=\sum_{i=0}^\infty{-1/2 \choose i}(-1)^i\rho^{2i}=\sum_{i=0}^\infty\frac{(2i-1)!!}{2^i i!}\rho^{2i}=\sum_{i=0}^\infty\frac{(2i-1)!!}{(2i)!!}\rho^{2i}, \]
and it remains to raise both sides to power $n$.

For equality, we must have an equality at each application of Prekopa-Leindler, which is characterized in \cite{dubuc}. Namely for each $I$, one has a.e. equalities $f=\lambda h(\bullet + C)$, $g=\lambda^{-1}h(\bullet - C)$. Already for a single multi-index $I$ this forces the existence of a diagonal matrix $A$ and $c\in\R$ such that a.e. one has $$\phi(x)=\frac 12 \langle Ax, x\rangle +c,\quad \psi(\xi)=\frac 12 \langle A^{-1}x, x\rangle -c.$$
\endproof

By the previous discussion, we immediately get
\begin{Corollary}
	For all unconditional Borel functions $f:\R^n\to \R_+$ it holds that \[ \widetilde M_\rho(f,0)\leq \frac{1}{n!}(2\pi)^n\rho^n(1-\rho^2)^{-\frac n 2}, \]
	with equality attained a.e.-uniquely by $f(x)=\lambda e^{-\frac12 \langle Ax, x\rangle}$ with some diagonal matrix $A>0$ and scalar $\lambda>0$. 
\end{Corollary}
\proof
Follows directly from Theorem \ref{thm:functional_unconditional}, since $\widetilde M_\rho(f,0)=\widetilde M_\rho(f)$ by unconditionality and the second part of Lemma \ref{lem:functional_invariance}.
\endproof

In particular, we proved a functional inequality of independent interest. 
\begin{Corollary}\label{cor:even_moments}
	Consider an unconditional Borel function $f=e^{-\phi}:\R^n\to \R_+$. Then $$ I_{2j}(f):=\int_{\R^n\times\R^n}\langle x,\xi\rangle ^{2j}e^{-\phi(x)-\mathcal L\phi(\xi)}dx d\xi\leq (2\pi)^n \frac{(n-2+2j)!!}{(n-2)!!(2j-1)!!}.$$
	Equality is attained a.e.-uniquely by multiples of gaussians $f(x)=\lambda e^{-\frac 12 \langle Ax, x\rangle}$, for some diagonal $A>0$ and scalar $\lambda>0$.
\end{Corollary}
The case of $j=0$ is the Blaschke-Santal\'o inequality for unconditional functions. The case of $j=1$ was established previously in \cite{huang_li}.
\proof
Denote $\psi=\mathcal L\phi$. In the proof of Theorem \ref{thm:functional_unconditional}, we have seen that $I_{2j}(e^{-\phi})$ is maximized a.e.-uniquely by multiples of gaussians. The maximal value $I_{2j}(G)$ can be computed directly, or alternatively it can be deduced from Theorem \ref{thm:functional_unconditional} as follows. We have $$\sum_{j=0}^\infty\frac{\rho^{2j}}{(2j)!}I_{2j}(G)= (2\pi)^n(1-\rho^2)^{-n/2},$$
and by Newton's binomial 
\[ (1-\rho^2)^{-\frac n2}=\sum_{j=0}^\infty{-n/2 \choose j}(-1)^j\rho^{2j}=\sum_{i=0}^\infty\frac{(n+2j-2)!!}{(n-2)!!(2j)!! }\rho^{2j}.\]
Consequently, 
$$I_{2j}(G)=(2\pi)^n(2j)!\frac{(n+2j-2)!!}{(n-2)!!(2j)!! }=(2\pi)^n\frac{(n-2+2j)!!}{(n-2)!!(2j-1)!!}.$$

\endproof

We can now deduce the corresponding inequalities for unconditional convex bodies. 
Define $$I_{2j}(K):=\int_{K\times K^o} \langle x,\xi\rangle^{2j}dxd\xi.$$
First, we prove two simple relations.

\begin{Lemma} For any convex body $K\in\mathcal K_0(\R^n)$,
	\begin{equation}\label{eq:gaussian_body} I_{2j}(K)=\frac{1}{((n+2j)!!)^2} \left(\frac 2\pi \right)^{\frac12(1-(-1)^n)}I_{2j}(e^{-\frac12 \|x\|_K^2}).\end{equation}
	If $K=-K$ then  
	\begin{equation}\label{eq:mahler_to_moments}\widetilde M_\rho(K)= \sum_{j=0}^\infty {n+2j \choose n}I_{2j}(K)\rho^{n+2j}.\end{equation} 
\end{Lemma}
\proof
We will use the change of variables $x=tp$, $\xi=s\eta$, with $p\in\partial K, \eta\in \partial K^o$. We have $dx=r^{n-1}drd\mu_K(p), d\xi =s^{n-1}dsd\mu_{K^o}(\eta)$.

We may write

\begin{align*}I_{2j}(e^{-\frac12 \|x\|_K^2})&=\int_{\R^n\times \R^n}\langle x, \xi\rangle ^{2j}e^{-\frac{\|x\|^2_K}{2}}e^{-\frac{\|\xi\|^2_{K^o}}{2}}dxd\xi\\&=\int_{[0,\infty)^2}e^{-\frac {r^2}{2}-\frac{s^2}{2}}r^{n+2j-1}s^{n+2j-1}drds\int_{\partial K\times\partial K^o}\langle p, \eta\rangle ^{2j}d\mu_K(p)d\mu_{K^o}(\eta)
	\\&=\left(\int_0^\infty r^{n+2j-1} e^{-r^2/2}dr\right)^2\int_{\partial K\times\partial K^o}\langle p, \eta\rangle ^{2j}d\mu_K(p)d\mu_{K^o}(\eta)\\
	&=2^{n+2j-2}\Gamma(\frac{n}{2}+j)^2\int_{\partial K\times\partial K^o}\langle p, \eta\rangle ^{2j}d\mu_K(p)d\mu_{K^o}(\eta).\end{align*}
On the other hand, 
\begin{align*}
	I_{2j}(K)&=\int_{K\times K^o}\langle x, \xi\rangle^{2j}dxd\xi\\&=\int_{\partial K\times \partial K^o} \langle p, \eta\rangle ^{2j}d\mu_K(p)d\mu_{K^o}(\eta)\int_{[0,1]^2}r^{n+2j-1}s^{n+2j-1}drds\\&=\frac{1}{(n+2j)^2}\int_{\partial K\times \partial K^o} \langle p, \eta\rangle ^{2j}d\mu_K(p)d\mu_{K^o}(\eta),
	\end{align*}
and it remains to combine the two equalities to find
\begin{align*} I_{2j}(K)&=\frac{1}{2^{n+2j-2}(n+2j)^2\Gamma(\frac n 2+j )^2}I_{2j}(e^{-\frac12\|x\|^2_K})\\&=\frac{1}{((n+2j)!!)^2} \left(\frac 2\pi \right)^{\frac12(1-(-1)^n)}I_{2j}(e^{-\frac12 \|x\|_K^2}).\end{align*}

Now assume $K=-K$. By Lemmas \ref{lem:body_vs_function} and \ref{lem:mahler_convex}, 
\begin{align*}
	\widetilde M_\rho(K)=\widetilde M_\rho(K,0)=\widetilde M_\rho(\mathbbm 1_K, 0)=&\frac{\rho^n}{n!}\int_{K\times\R^n}e^{-\|\xi\|_{K^o}+\rho\langle x,\xi\rangle}dxd\xi\\&=\frac{\rho^n}{n!}\int_{K\times\R^n}e^{-\|\xi\|_{K^o}}\cosh(\rho\langle x,\xi\rangle)dxd\xi
	\\&=\sum_{j=0}^\infty \frac{1}{n!(2j)!}\rho^{n+2j}\int_{K\times\R^n}\langle x, \xi\rangle^{2j}e^{-\|\xi\|_{K^o}}dxd\xi.
\end{align*} 
We may write 
\begin{align*}I_{2j}(K)=\int_{K\times K^o}\langle x,\xi\rangle^{2j}dxd\xi&=\int_{K\times\partial K^o}\langle x,\eta\rangle^{2j}dxd\mu_{K^o}(\eta)\int_0^1s^{2j+n-1}ds\\&=\frac{1}{n+2j}\int_{K\times\partial K^o}\langle x,\eta\rangle^{2j}dxd\mu_{K^o}(\eta).\end{align*}
Therefore, 
\begin{align*}\int_{K\times\R^n}\langle x, \xi\rangle^{2j}e^{-\|\xi\|_{K^o}}dxd\xi &=\int_{K\times [0,\infty)\times \partial K^o} s^{2j+n-1}e^{-s}\langle x,\eta\rangle ^{2j}dxds d\mu_{K^o}(\eta)\\&=(2j+n-1)!\int_{K\times\partial K^o}\langle x, \eta\rangle^{2j}dxd\mu_{K^o}(\eta)\\&=(2j+n)!I_{2j}(K),\end{align*}
and so
\[	\widetilde M_\rho(K) = \sum_{j=0}^\infty \frac{(n+2j)!}{n!(2j)!}\rho^{n+2j}I_{2j}(K). \]
\endproof

\begin{Corollary}\label{cor:bodies_moments}
	For $K\subset \R^n $ convex and unconditional, 
	$$ I_{2j}(K):=\int_{K\times K^o} \langle x,\xi\rangle^{2j}dxd\xi\leq \frac{(2\pi)^n}{n+2j}\frac{1}{ (n+2j)!!(n-2)!!(2j-1)!!} \left(\frac2\pi\right)^{\frac{1-(-1)^n}{2} }.$$ Equality is attained uniquely by ellipsoids.
\end{Corollary}
\proof
We apply Corollary \ref{cor:even_moments} to $f(x)=e^{-\phi(x)}=e^{-\frac12 \|x\|_K^2}$.  By eq. \eqref{eq:gaussian_body}, $I_{2j}(K)$ is uniquely maximized by ellipsoids, and the maximal value is readily computed.
\endproof

\begin{Corollary}\label{cor:bodies_mahler}
	For any $0\leq  r<\infty$, and $K\subset \R^n$ convex and unconditional, $M_r(K)\leq  M_r(B^n)$. Equality is attained uniquely by ellipsoids.
\end{Corollary}
\proof
We only need to note that by Corollary \ref{cor:bodies_moments}, each summand of \eqref{eq:mahler_to_moments} is uniquely maximized by ellipsoids.
\endproof

\textit{Proof of Corollary \ref{cor:entropy}}. Since $d^H=\frac12(d^F+d^{RF})\geq \frac12d^F$, it follows that $B^H_r(q, K)\subset B^F_{2r}(q, K)$.  Write $c_n$ for a constant depending only on $n$ which may change between occurrences. By Proposition \ref{prop:centro_affine} we have $\vol^F(B^F_R(q,B^n))\leq c_ne^{\frac{n-1}{2}R }$. Using the Rogers-Shephard inequality as in the proof of Corollary \ref{cor:hilbert_volume_sec3}, and applying Corollary \ref{cor:bodies_mahler}, we find
\begin{align*} \vol_K^H(B_r^H(q,K))&\leq c_n \vol_K^F(B_r^H(q,K))\leq c_n\vol_K^F(B_{2r}^F(q,K))\leq  c_n\vol_K^F(B_{2r}^F(q,B^n))\\&\leq c_ne^{\frac{n-1}{2}2r}=c_ne^{(n-1)r},\end{align*}
completing the proof. \qed

Using the same methods as before, we can also show the following.

\begin{Theorem}
	
	Let $f=e^{-\phi}:\R^n_+\to \R_+$ be a Borel function. Then
	\begin{itemize}
		\item For all $j\geq 0$, $$I^+_j(f):=\int_{\R_+^n\times\R_+^n}\langle x,\xi\rangle^{j} e^{-\phi(x)-\mathcal L\phi(\xi)}dx d\xi\leq I^+_j(e^{-\frac{|x|^2}{2}}).$$ 
		\item If $n=1$,  
		$$I_{2j+1}^+(f)=\int_0^\infty x^{2j+1} e^{-\phi(x)}dx  \int_0^\infty \xi^{2j+1}e^{-\mathcal L\phi(\xi)} d\xi\leq (2j)!!^2.$$ 
		\item For all $0\leq \rho<1$, $$\int_{\R_+\times\R_+}e^{-\phi(x)-\mathcal L\phi (\xi)}\sinh (\rho x\xi)dxd\xi\leq\frac{1}{\sqrt{1-\rho^2}}\arctan\frac{\rho}{\sqrt{1-\rho^2}}.$$
	\end{itemize}
	
	Equality in all cases is a.e.-uniquely attained by $\lambda e^{-\frac12 \langle Ax, x\rangle}$, where $A>0$ is diagonal and $\lambda>0$.  
	
\end{Theorem}
\proof
The first two inequalities can be obtained by trivially adjusting the proof of Theorem \ref{thm:functional_unconditional}. For the last inequality, we should show that
$$\sum_{j=0}^\infty\frac{(2j)!!^2}{(2j+1)!}\rho^{2j+1}=\sum_{j=0}^\infty\frac{(2j)!!}{(2j+1)!!}\rho^{2j+1}=\frac{1}{\sqrt{1-\rho^2}}\arctan\frac{\rho}{\sqrt{1-\rho^2}}.$$ This can be done by observing  that
$$ \frac{(2j)!!}{(2j+1)!!}=\int_0^{\pi/2}\sin^{2j}tdt,$$
and so \begin{align*}\sum_{j=0}^\infty\frac{(2j)!!}{(2j+1)!!}\rho^{2j+1}&=\int_{0}^{\frac \pi 2}(\rho\sin t)^{2j+1}dt=\int_0^{\frac \pi 2}\frac{\rho \sin t}{1-\rho^2\sin^2t}dt\\&=\int_0^1 \frac{2\rho ds}{s^2+2(1-2\rho^2)s +1}\\&=\frac{1}{\sqrt{1-\rho^2}}\left(\arctan \frac{\sqrt{1-\rho^2}}{\rho}-\arctan \frac{1-2\rho^2}{2\rho \sqrt{1-\rho^2}}\right)\\&=\frac{1}{\sqrt{1-\rho^2}}\arctan \frac{\rho}{\sqrt{1-\rho^2}}.\end{align*}
Equality cases are settled using \cite{dubuc} as in Theorem \ref{thm:functional_unconditional}.
\endproof
As a sidenote, the discussion above suggests an extension of the centro-affine surface area to functions.

\begin{Definition}
	Let $f:\R^n\to[0,\infty)$ be a Borel function.  Its centro-affine surface area is
	\begin{align*}\Omega_n(f)&:=\frac{1}{n!}\limsup_{\rho\to1^-} (1-\rho^2)^{n/2}\displaystyle\int_{\R^n\times \R^n} e^{-\phi(x)-\mathcal L \phi (\xi)+\rho\langle x,\xi\rangle}dxd\xi.\end{align*}
\end{Definition}
Evidently, $\Omega_n(f)=\Omega_n(\lambda f\circ A)$ for any $\lambda> 0$ and $A\in\GL_n(\R)$.

\section{Regularizing the total Funk volume} \label{sec:projective_total_volume}
We use the notation of subsection \ref{sec:volume_funk}. Recall that we use the standard Euclidean structure on $V=\R^{n+1}$, and identify the double cover of $\PP$ and $\PP^\vee$ with the unit sphere. Let $x_0\in x$ be a unit vector, well-defined up to sign. We write $\sigma_x, \sigma_\xi$ for the standard Euclidean measure on the corresponding copy of $S^2$.

\begin{Lemma}\label{lem:lemma91}For $g\in \PGL(V)$, set $\chi_g(x,\xi)=|gx_0||g^{-*}\xi_0|$.  Then
	\[g_*(\sigma_x\sigma_\xi)=\chi_{g^{-1}}^{-n-1}\sigma_x\sigma_\xi\]
\end{Lemma} 
\proof
Denote by $\Stab^+(x,\xi)\subset\Stab(x,\xi)$ the subgroup of the stabilzer in $\PGL(V)$ that has $\sign \det g|_x=\sign \det g^*|_\xi$. Then $\Stab^+(x,\xi)$-equivariantly we have $$\Dens(T_{x,\xi}(\PP\times \PP^\vee))=\Dens^*(x)^{n+1}\otimes\Dens^*(\xi)^{n+1}=x^{n+1}\otimes \xi^{n+1},$$
and so
\[g_*(\sigma_x\sigma_\xi)(gx,g^{-*}\xi)=(|gx_0||g^{-*}\xi_0|)^{n+1}\sigma_x\xi=\chi_g(x,\xi)^{n+1}\sigma_x\sigma_\xi.\]
Hence
\[g_*(\sigma_x\sigma_\xi)(x,\xi)=\chi_g(g^{-1}x,g^*\xi)^{n+1}\sigma_x\sigma_\xi=\chi_{g^{-1}}^{-n-1}(x,\xi)\sigma_x\sigma_\xi,\]
as claimed.
\endproof
Define $f_0\in C(\PP\times\PP^\vee)$ by $f_0(x,\xi)=|\langle x_0,\xi_0\rangle|$.
It is easy to check that \[ g^*f_0=\chi_g^{-1}f_0.\]
Define $\mu_z(x,\xi)=f_0(x,\xi)^z\sigma_x\sigma_\xi\in \mathcal M(\PP\times\PP^\vee\setminus \mathcal Z)$. For $\textrm{Re}z>-1$, it extends as an absolutely continuous measure on $\PP\times\PP^\vee$.
\begin{Lemma}\label{lem:measure_push}
	$g_*\mu_z=\chi_{g^{-1}}^{-(n+1)-z}\mu_z$.
\end{Lemma}
\proof
It is easy to check that $g_*\mu_z=(g^{-1})^*f_0^z\cdot g_*(\sigma_x\sigma_\xi)$. It remains to put the previous computations together.
\endproof
As expected from Proposition \ref{prop:invariant_measure}, we see that $z=-(n+1)$ corresponds to the invariant measure $\mu$.

Henceforth we use the same notation for a point on $S^n$ and the corresponding line in $\R\PP^n$. The action of $\GL(V)$ on $S^n$ is $g(x)=\frac{gx}{|gx|}$.

Let $K\subset S^n\subset V=\R^{n+1}$ be a proper convex compact set, and $K^\vee\subset S^n$ its polar body, namely $K^\vee=\{\xi\in S^{n}: \forall x\in K, \langle x,\xi\rangle\geq 0\}$. They correspond to a dual pair of convex sets in $\PP$ and $\PP^\vee$, also denoted $K, K^\vee$.

Define the \emph{Beta function} of $K\subset S^n$ by 
\[\C\ni z\mapsto B_K(z)=\int_{K\times K^\vee}\langle x,\xi\rangle ^z\sigma_x\sigma_\xi =\langle  \mathbbm1_{K\times K^\vee}, \mu_z\rangle. \]
Clearly, $B_K(z)=B_{K^\vee}(z)$. Assume henceforth that $K$ is smooth and strictly convex. 
\begin{Lemma}\label{lem:convergence_domain}
	The integral $B_K(z)$ is convergent if and only if $\mathrm{Re} z>-\frac{n+3}{2}$.
\end{Lemma}
\proof
The zeros of $\langle x,\xi\rangle $ on $K\times K^\vee$ are $\mathcal Z\cap K\times K^\vee$, which is an $(n-1)$-dimensional submanifold in $\partial K\times\partial K^\vee$. Let $x_1,\dots,x_{n-1}$ be local coordinates on $\partial K$, and $\xi_1,\dots, \xi_{n-1}$ the Legendre-dual coordinates on $\partial K^\vee$, that is $\mathcal Z\cap (K\times K^\vee)=\{\xi_j=x_j\}$. Let $x_n$ be the remaining coordinate near $\partial K$, so that locally $ K=\{x_n\geq 0\}$, and $K^\vee=\{\xi_n\geq 0\}$.
We put $\xi=x+y$, and use the notation $v^j=(v_i)_{i=1}^{j}$. When restricted to $\{x_n=\xi_n=0\}$, $f_0( (x^{n-1},0), (\xi^{n-1},0)) =\langle A_x y^{n-1}, y^{n-1}\rangle+O(|y^{n-1}|^3)$, for some $A_x>0$. It follows that in general, $$f_0(x,\xi)=\langle A_{x}y^{n-1}, y^{n-1}\rangle+b_xx_n+c_x\xi_n +O(|y^{n-1}|^3+|y^{n-1}|(x_n+\xi_n)+x_n^2+\xi_n^2),$$ where $A_x, b_x, c_x>0$ are uniformly bounded from above and below by strict convexity.

Thus $\int_{K\times K^\vee} f_0(x,\xi)^z$ converges if and only if

\[\int_{y_1^2+\dots+y_{n-1}^2<1}\int_{x_n, \xi_{n}=0}^1 ( y_1^2+\dots+y_{n-1}^2+x_n+\xi_{n})^z<\infty,  \]
which immediately reduces to $2(\mathrm{Re}z+2)+(n-2)>-1\iff \Re z>-\frac{n+3}{2}$.
\endproof

\emph{In all  the following we assume $n=2$}, as already this case appears to be interesting and involved. We proceed to prove Theorem \ref{main:beta}.

Denote $D_+=\{(x,\xi); \langle x,\xi\rangle> 0\}\subset S^2\times S^2$. The group $\SO(3)$ acts diagonally on $S^2\times S^2$, leaving $D_+$ invariant. By $\Omega^{i,j}(S^2\times S^2)$ we denote the differential forms of bidegree $(i,j)$. We will use the function $R_z(t)=\frac{1-t^{z}}{1-t^2}$, which is smooth in $t>0$ and analytic in $z\in\mathbb C$.

\begin{Proposition}\label{prop:primitive} For $z\in \C$, $z\neq -1$, define $\omega_z\in{ \Omega^{1,1}(D_+)^{\SO(3)}  }$  by \[\omega_z:=\frac{1}{z+1}R_{z+1}(\langle x,\xi\rangle)\langle dx,\xi\rangle\wedge \langle x, d\xi\rangle.\]
	Then $d_xd_\xi\omega_z=\langle x,\xi\rangle^z \sigma_x\sigma_\xi$.
\end{Proposition}
\proof
Rather than verify the claimed equality, we will construct from scratch the form $\omega_z$ with the required derivative. We write $S^2_x$, $S^2_\xi$ for the corresponding copies of the sphere.

Fix $\xi$ and use it as the North Pole with associated polar coordinates $\theta, \phi$ of $x\in S^2$, where $\cos\theta=\langle x,\xi\rangle$, and $\phi$ is the azimuth. We have $\sigma_x=\sin\theta d\theta d\phi$. It holds that $d\phi(u)=-\frac{1}{\sin\theta}d\theta(x\times u)$, and differentiating $\cos\theta=\langle x,\xi\rangle$ we find
\begin{equation}\label{eq:dphi} -\sin \theta d\theta=\langle dx, \xi\rangle\Rightarrow d\phi=\frac{1}{\sin^2\theta}\langle x\times dx, \xi\rangle=\frac{\langle x\times dx, \xi\rangle}{1-\langle x,\xi\rangle^2} \end{equation}
The forms $d\theta$, $\langle x\times dx,\xi\rangle=\sin^2\theta d\phi$ are linearly independent in $S^2\setminus \{\pm \xi\}$, and vanish at $\pm\xi$. 
Assume a 1-form $\zeta_z$ on $S^2_x$ of the form $\zeta_z=b_z(\theta)d\phi$ is given, such that $d_\xi\omega_z= \zeta_z \sigma_\xi.$  Then
\[ d\zeta_z=b_z'(\theta) d\theta d\phi+ c_0 b_z(0)\delta_\xi,  \]
where $c_0=2\pi$ and $\delta_\xi$ is the delta measure at $\xi$.
Thus if $b_z(0)=0$, and $b_z'(\theta)=\sin\theta\cos^z\theta$, the desired equation would be satisfied. For this we take 
\[b_z(\theta)=\frac{1-\cos^{z+1}\theta}{z+1	} \]

Now fix $x$. We should find an  $\R^3$-valued $\Stab(x)=\SO(2)$-equivariant $1$-form $\eta_z$ on $S^2_\xi$ with 
\begin{equation}\label{eqn:eta_differential}d_\xi\eta_z=\frac{1-\langle x,\xi\rangle^{z+1}}{1-\langle x,\xi\rangle^2}  \sigma_\xi \xi, \end{equation}
whence we may take $\omega_z= \frac{1}{z+1}  \langle x\times dx, \eta_z\rangle $.
It is easy to verify that $$\eta_z=a_z(\theta)d\phi \cdot x+c_z(\theta)d\theta \cdot x\times\xi,$$ where $c_z(\theta)=\frac{1-\cos^{z+1}\theta}{\sin\theta}$ and $a_z(\theta)=\int_0^\theta c_z(\tau)\cos\tau d\tau$ satisfies \eqref{eqn:eta_differential}. The verification boils down to the equality 
$$\frac{1-\cos^{z+1}\theta}{\sin\theta}\cos\theta d\theta d\phi\cdot x-\frac{1-\cos^{z+1}\theta}{\sin\theta}d\theta \cdot x\times d\xi=\frac{1-\cos^{z+1}\theta}{\sin\theta}d\theta d\phi\cdot \xi,$$ which in turn follows from \eqref{eq:dphi} applied to $S^2_\xi$, by taking scalar products with $x,\xi$ and $x\times\xi$.

We find that \[\omega_z=-\frac{1}{z+1}\langle x\times dx, \eta_z\rangle=\frac{1}{z+1}\frac{1-\langle x,\xi\rangle^{z+1}}{1-\langle x,\xi\rangle^2}\langle dx,\xi\rangle\wedge \langle x,  d\xi\rangle.\]

\endproof

It follows by the Stokes theorem and using continuity to pass to the boundary of $D_+$, that for $\Re z\geq 0$,
\begin{equation}\label{eq:stokes_beta}B_K(z) = \int_{ K\times K^\vee} \langle x,\xi\rangle^z = \int_{\partial K\times \partial K^\vee} \omega_z. \end{equation}

\begin{exmp}  Let us compute $B_K(z)$ for the spherical disc $K=\{x_3\geq \frac{1}{\sqrt 2}\}$, satisfying $K=K^\vee$. 
Parametrize $\partial K$ by $x(\alpha)=\frac{1}{\sqrt 2}(\cos \alpha, \sin\alpha, 1)$, and $\partial K^o$ by $\xi(\beta)=\mathcal L(x(\beta))=\frac{1}{\sqrt 2}(-\cos \beta, -\sin\beta, 1)$, $0\leq \alpha,\beta\leq 2\pi$.
Then: \[\langle x,\xi\rangle =\frac12(1-\cos(\alpha-\beta))=\sin^2\frac{\alpha-\beta}{2}\]
\[ dx=\frac{1}{\sqrt 2}(-\sin\alpha,\cos\alpha,0)d\alpha,\quad  d\xi=\frac{1}{\sqrt 2}(\sin\beta,-\cos\beta,0)d\beta \]
\[\langle dx,\xi\rangle = \frac12\sin(\alpha-\beta)d\alpha,\quad \langle d\xi,x\rangle =\frac12\sin(\beta-\alpha)d\beta  \]
Thus
\begin{align*} B_K(z)&=-\frac{1}{4(z+1)} \int_0^{2\pi}\int_0^{2\pi}\frac{1-(\sin^2\frac{\alpha-\beta}{2})^{z+1}}{1-\sin^4\frac{\alpha-\beta}{2}}\sin^2(\alpha-\beta)d\alpha d\beta 
\\& =-\frac{1}{2(z+1)}\int_{-\pi}^{\pi}(2\pi-2|\gamma|)\frac{1-(\sin^2\gamma)^{z+1}}{1-\sin^4\gamma}\sin^22\gamma d\gamma\\&=
-\frac{8}{z+1}\int_0^{\pi}(\pi-\gamma)\frac{1-(\sin^2\gamma)^{z+1}}{1+\sin^2\gamma}\sin^2\gamma d\gamma\\&=
-\frac{8}{z+1}\int_0^{\pi}\frac \pi 2\frac{1-(\sin^2\gamma)^{z+1}}{1+\sin^2\gamma}\sin^2\gamma d\gamma
\\&=
-\frac{8\pi}{z+1}\int_0^{\frac\pi2}\frac{1-(\sin^2\gamma)^{z+1}}{1+\sin^2\gamma}\sin^2\gamma d\gamma, \end{align*}
where the penultimate equality is by symmetrizing the integrand around $\frac\pi 2$. 

Putting $t=\sin^2\gamma$ we get
\begin{align*} B_K(z)&=-\frac{4\pi}{z+1}\int_0^1\frac{1-t^{z+1}}{1+t}\frac{\sqrt t}{\sqrt{1-t}}dt=-\frac{4\pi}{z+1}(\pi(1-\frac{1}{\sqrt 2})-\int_0^1 \frac{t^{z+\frac 32}dt}{(1+t)\sqrt{1-t}} ) 
\\&=-\frac{4\pi}{z+1}\left(\pi(1-\frac{1}{\sqrt 2})-\sqrt\pi \frac{\Gamma(z+\frac 52)}{\Gamma(z+3)}{}_2F_1(1,z+\frac52,z+3;-1)\right) 
.\end{align*}
In particular, $B_K(-3)=2\pi^2$.
\end{exmp}

\begin{Proposition}
	For $K\subset S^2$ smooth and strictly convex, $B_K(z)$ extends as a meromorphic function with simple poles that are contained in $\{-\frac52, -\frac72,-\frac 92,\dots\}$. 
\end{Proposition}
\proof
Let $\mathcal L:\partial K\to \partial K^\vee$ be the Legendre transform, defined by  $\mathcal L x=\partial K^\vee\cap x^\perp$. It points in the unique direction perpendicular to both $x$ and $T_x\partial K$: since $\langle x,\xi\rangle:\partial K\times\partial K^\vee\to \R$ is smooth and non-negative, its zeros are all minima and hence critical points. That is, for fixed $\xi=\mathcal L x$, we have $\langle dx,\xi\rangle =0$. Parametrize $\partial K$ by arc-length, denoted $x(\alpha)$, and $\partial K^\vee$ by $\xi(\beta)=\mathcal L(x(\beta))$. It then holds that 
\[ \langle x(\alpha), \xi(\beta)\rangle =(\beta-\alpha)^2 u(\alpha,\beta)\]
with $u$ smooth, and by strict convexity $u>0$. Similarly we have $$\langle x'(\alpha),\xi(\beta)\rangle =(\beta-\alpha) v_1(\alpha,\beta),\quad \langle \xi'(\beta),x(\alpha)\rangle =(\beta-\alpha) v_2(\alpha,\beta),$$ with $v_1,v_2$ smooth and non-zero when  $|\alpha-\beta|\leq \epsilon$. Here $\epsilon>0$  is fixed, and depends on the lower bound on the curvatures of $K, K^\vee$. 

Recall eq. \eqref{eq:stokes_beta}: 
\[B_K(z) =-\frac{1}{z+1}\int_{\partial K\times \partial K^\vee} \frac{1-\langle x,\xi\rangle^{z+1}}{1-\langle x,\xi\rangle^2}\langle dx,\xi\rangle\wedge \langle d\xi,x\rangle .\]
Writing $\beta=\alpha+\phi$, it follows that up to an addition of an entire function of $z$ and the factor $\frac{1}{z+1}$, $B_K(z)$ equals

\begin{align*}B^\epsilon_K(z)&:=\int_{-\epsilon}^\epsilon d\phi\int_{\alpha=0}^{2\pi}\frac{1-\langle x(\alpha),\xi(\alpha+\phi)\rangle^{z+1}}{1-\langle x(\alpha),\xi(\alpha+\phi)\rangle^{2}} \langle x'(\alpha), \xi(\alpha+\phi)\rangle
	\langle \xi'(\alpha+\phi), x(\alpha)\rangle  d\alpha 
	\\ &=h_\epsilon(K)-\int_{-\epsilon}^\epsilon|\phi|^{2z+4} d\phi\int_{\alpha=0}^{2\pi}  u(\alpha,\alpha+\phi)^{z+1}\tilde v(\alpha,\phi)d\alpha ,
\end{align*}
where $h_\epsilon(K)$ is a constant depending only on $K$, and $ \tilde v$ is a smooth non-zero function.

The inner integral, denoted $f_z(\phi)$, is an analytic family of smooth functions as $z\in\mathbb C$. The function $\int_{-\epsilon}^\epsilon|\phi|^{2z+4} f_z(\phi)d\phi$ is then meromorphic in $\C$, with simple poles that are contained in a subset of $\{-\frac52,-\frac 72,\dots\}$, corresponding to the simple poles of the meromorphic family  $\{|\phi|^{2z+4}\}$ of even homogeneous distributions on $\R$. Also by Lemma \ref{lem:convergence_domain}, $z=-1$ is an analytic point of $B_K(z)$. This completes the proof of the meromorphic extendibility of $B_K(z)$ with poles as stated.
\endproof

\begin{Proposition}
	The value $B_K(-3)$ is a projective invariant of $K\subset\R\PP^2$.
\end{Proposition}
\proof
Recall the symbol $\chi_g(x,\xi)$ defined in Lemma \ref{lem:lemma91}. For any fixed $g\in \PGL(V)$, we have by Lemma \ref{lem:measure_push} the following equality for $\textrm{Re}z>-1$:
\begin{align}\label{eq:invariance}&\langle \mathbbm 1_{g^{-1}K\times (g^{-1}K)^\vee}, \mu_z\rangle=\langle g^*\mathbbm1_{K\times K^\vee},\mu_z\rangle=\langle \mathbbm1_{K\times K^\vee}, g_*\mu_z\rangle\\=&\langle \chi_{g^{-1}}^{-3-z}\mathbbm 1_{K\times K^\vee}, \mu_z\rangle.\nonumber\end{align}
We already established that one end (and hence both) of this equation admits a meromorphic extension in $z\in \mathbb C$, which is analytic at $z=-3$. Let us verify that the value of the right hand side at $z=-3$ equals $B_K(-3)$, implying projective invariance.

 Fix $g\in \PGL(3)$, and consider $h(z)=\langle \chi_{g^{-1}}^{-3-z} \mathbbm 1_{K\times K^\vee}, \mu_z\rangle=\int_{K\times K^\vee} \chi_{g^{-1}}^{-3-z} d\mu_z$. Write $\psi(x,\xi)=\chi_{g^{-1}}(x,\xi)$, and note $\psi$ is a smooth positive function. We now use Proposition \ref{prop:primitive} and the Stokes theorem to write 

\begin{align*}h(z)&= \int_{K\times K^\vee} \psi(x,\xi)^{-z-3}d_xd_\xi\omega_z\\&=\int\limits_{\partial K\times K^\vee}\psi^{-z-3}d_\xi\omega_z+(z+3)\int\limits_{K\times K^\vee}\psi^{-z-4} d_x \psi\wedge d_\xi \omega_z 
	\\ &=\int\limits_{\partial K\times \partial K^\vee}\psi^{-z-3}\omega_z- (z+3)\int\limits_{K\times K^\vee} d_\xi(\psi^{-z-4}d_x \psi)\wedge\omega_z \\ &+(z+3)\int\limits_{\partial K\times K^\vee}\psi^{-z-4}d_\xi \psi\wedge\omega_z +(z+3)\int\limits_{ K\times \partial K^\vee}\psi^{-z-4}  d_x \psi\wedge \omega_z\end{align*}

The first summand admits a meromorphic extension by the same proof as that for $B_K$, and its value at $z=-3$ is  $\int_{\partial K\times \partial K^\vee}\omega_{-3}=B_K(-3)$. Consider next the second summand; by Lemma \ref{lem:convergence_domain}, the form $\omega_z$ is $L^1$ in $K\times K^\vee$ for $\Re z>-\frac72$, implying the second summand is analytic in this range, and vanishes at $z=-3$. 

The last two summands are treated similarly, so let us focus on the first one. Written explicitly, it is

$$\frac{z+3}{z+1}\int_{\partial K\times K^\vee} \frac{1-\langle x,\xi\rangle^{z+1}}{1-\langle x,\xi\rangle^2}\langle dx,\xi\rangle\wedge\langle x,d\xi\rangle\wedge  \psi^{-z-4}d_\xi\psi.$$
It is easy to see that up to a summand that is analytic near, and vanishes at, $z=-3$, the integral equals
$$I=\frac{z+3}{z+1}\int_{\partial K\times K^\vee} \langle x,\xi\rangle^{z+1}\langle x,d\xi\rangle \wedge\langle dx,\xi\rangle\wedge \psi^{-z-4}d_\xi\psi.$$
It holds that $$\frac{1}{z+2}d(\langle x,\xi\rangle^{z+2})\wedge\langle dx,\xi\rangle\wedge \psi^{-z-4}d_\xi\psi=\langle x,\xi\rangle^{z+1}\langle x, d\xi\rangle \wedge \langle dx, \xi\rangle\wedge \psi^{-z-4}d_\xi\psi.$$
Consequently we can integrate by parts to obtain 
\begin{align*}I=I_1+I_2&=\frac{z+3}{(z+1)(z+2)}\int\limits_{\partial K\times\partial K^\vee}\langle x,\xi\rangle ^{z+2} \langle dx,\xi\rangle\wedge\psi^{-z-4}d_\xi\psi\\&-\frac{z+3}{(z+1)(z+2)}\int_{\partial K\times K^\vee} \langle x,\xi\rangle ^{z+2}  d(\langle dx,\xi\rangle \wedge\psi^{-z-4}d_\xi\psi). \end{align*}

Let us rewrite the first summand as
\begin{align}\nonumber I_1&=\frac{1}{(z+1)(z+2)}\int\limits_{\partial K\times \partial K^\vee}d_x(\langle x,\xi\rangle^{z+3})\wedge \psi^{-z-4}d_\xi\psi\\&=- \frac{1}{(z+1)(z+2)}\int\limits_{\partial K\times\partial K^\vee} \langle x,\xi\rangle^{z+3}d_x(\psi^{-z-4}d_\xi\psi)\label{eq:I1_2}  \end{align}

We use the arc-length parametrization $x(\alpha)$, $0\leq \alpha\leq L$ for $\partial K$, and $\xi (\beta)=\mathcal L(x(\beta))$, $0\leq \beta\leq L$ for $\partial K^\vee$. Recall that by strict convexity, $\langle x(\alpha), \xi(\beta)\rangle=|\beta-\alpha|^2u(\alpha, \beta)$ with $u>0$ smooth. We will use $C$ to denote a positive constant that only depends on $K$ and $g$, which may differ between occurrences. 

Consider the integral in \eqref{eq:I1_2}. As $z\to-3$, the integrand is majorized by $C|\langle x, \xi\rangle|^{-1/4}$, which is integrable on $\partial K\times\partial K^\vee$.  Thus by Lebesgue's dominated convergence theorem, $$\lim _{z\to -3} I_1=-\frac12\int_{\partial K\times\partial K^\vee}d_x(\psi^{-1}d_\xi\psi) =0.$$

Next we consider $I_2$. Let $P: S^2\to \partial K^\vee$ be the nearest-point projection, well-defined in some open neighborhood $U$ of $\partial K^\vee$. Define $K^\vee_\epsilon=\{\xi\in U: \langle \xi, \mathcal L(P\xi)\rangle \leq \epsilon\}$, which is a one-sided tube around $\partial K^\vee$.
We parametrize $\partial K\times K^\vee_\epsilon$ by $(x(\alpha), \xi(\beta, s))$, $0\leq \alpha,\beta\leq L$, $0\leq s\leq \epsilon$, where $\xi=\xi(\beta, s)$ is the point that has $P(\xi)=\xi(\beta)$, $\langle \xi, \mathcal L (P(\xi))\rangle=s$. 

Write $\eta=d(\langle dx,\xi\rangle \wedge\psi^{-z-4}d_\xi\psi)$. We claim that the integral 
\begin{equation}\label{eq:integral}\int_{\partial K\times K^\vee}\langle x,\xi\rangle^{z+2}\eta\end{equation}
converges for $\Re z>-\frac 72$. 

It is enough to consider the integral \eqref{eq:integral} in the smaller domain $\{|\alpha-\beta|\leq \epsilon, 0\leq s\leq \epsilon\}$. There $\langle x(\alpha), \xi(\beta, s)\rangle \geq C\max((\beta-\alpha)^2, s)\geq C((\beta-\alpha)^2+s) $. 
It remains to check that for $\Re z>-\frac72$,
\[ \int_0^\epsilon\int_0^\epsilon\frac{dsd\phi}{(s+\phi^2)^{-z-2}}<\infty,  \]
which is clear. Consequently, $\lim_{z\to -3} I_2=0$ due to the $(z+3)$ factor. This concludes the proof of the statement.

\endproof

Let us derive some equivalent expressions for $B_K(-3)$, working towards an explicit geometric formula. One has

\begin{align*}B_K(-3)&=-\frac12\left. \int_{\partial K\times\partial K^\vee}\frac{1-\langle x,\xi\rangle^{z+1}}{1-\langle x,\xi\rangle ^2}\langle dx,\xi\rangle\wedge\langle x,d\xi\rangle\right|_{z=-3} \\&=\frac12 \left.\int_{\partial K\times\partial K^\vee}\langle x,\xi\rangle^{z+1} \frac{1-\langle x,\xi\rangle^{-z-1}}{1-\langle x,\xi\rangle ^2}\langle dx,\xi\rangle\wedge\langle x,d\xi\rangle\right|_{z=-3}\\&=\frac12  \left.\int_{\partial K\times\partial K^\vee}\langle x,\xi\rangle^{z+1}\langle dx,\xi\rangle\wedge\langle x,d\xi\rangle\right|_{z=-3}
	\\&=\frac1{2(z+2)}\left.\int_{\partial K\times\partial K^\vee}\langle dx,\xi\rangle \wedge d(\langle x, \xi\rangle^{z+2})\right|_{z=-3}
	\\&=\frac1{2}\left.\int_{\partial K\times\partial K^\vee}\langle x,\xi\rangle^{z+2}\langle dx,d\xi\rangle\right|_{z=-3}.
\end{align*}
\begin{align}B_K(-3)&=\frac1{2}\left.\int_{\partial K\times\partial K^\vee}\frac{\langle dx,d\xi\rangle}{\langle x,\xi\rangle^z}\right|_{z=1}
	\label{eq:mahler_volume} 
	\\&=\frac12\left.\int_{\partial K\times\partial K^\vee} \frac{\langle dx,\xi\rangle\wedge \langle x,d\xi\rangle}{\langle x,\xi\rangle^z}\right|_{z=2}\label{eq:mahler_volume2} 
\end{align}

\begin{Proposition}
	Let $\partial K$ be parametrized by arc-length as $x(\alpha)$. The geodesic curvature of $\partial K$ is $\kappa(\alpha)=\det (x(\alpha), x'(\alpha), x''(\alpha))$. Then
	$$B_K(-3)= \frac12 \lim_{\epsilon\to 0}\int_{x\in\partial K, \xi\in \partial K^\vee, \langle x,\xi\rangle \geq\epsilon^2}\frac{\langle dx,\xi\rangle\wedge\langle x,d\xi\rangle}{\langle x,\xi\rangle^2}+\frac{2\sqrt 2}{\epsilon}\int_{\partial K}\sqrt {\kappa(\alpha )}d\alpha.$$
\end{Proposition}
\proof
We carry out the regularization explicitly, using formula \eqref{eq:mahler_volume2}. Put $L=\textrm{Length}(x(\alpha))$, $\xi(\beta)=\mathcal L(x(\beta))=x(\beta)\times x'(\beta)$, and let $\alpha,\beta\in\R$ be $L$-periodic parameters.

We will write $x(\alpha)>\xi(\beta)$ if $\beta<\alpha<\beta+\frac{L}{2}$. Define for small values $|t|\leq \epsilon$
\[g_K(t)=\left\{\begin{array}{cc}\int_{\langle x,\xi\rangle =t^2, x>\xi}\langle dx,\xi\rangle&,\qquad t\geq0 \\ \int_{\langle x,\xi\rangle =t^2, x<\xi}\langle dx,\xi\rangle&,\qquad t<0 \end{array}\right.\]
Clearly $g_K$ is smooth in $[-\epsilon,\epsilon]$, and $g_K(0)=0$.

Let $S(x,\xi)=\sqrt{\langle x,\xi\rangle}$ be a signed square root defined on $\partial K\times \partial K^\vee$ , which is non-negative for $x\geq \mathcal L(\xi)$, and non-positive for $x\leq L(\xi)$. Put $t=S(x,\xi)$, so $t^2=f_0(x,\xi)=\langle x,\xi\rangle$ and $2tdt=\langle dx,\xi\rangle +\langle x,d\xi\rangle$, and $\langle dx,\xi\rangle \wedge\langle x,d\xi\rangle=-2tdt\langle dx,\xi\rangle$.

It then holds that
\begin{align*}B_K(-3)&= \frac12 \int_{\langle x,\xi\rangle \geq\epsilon^2}\frac{\langle dx,\xi\rangle\wedge\langle x,d\xi\rangle}{\langle x,\xi\rangle^2} -  \int_0^{\epsilon}\frac{g_K(t)-g_K(-t)-2g_K'(0)t}{t^3}dt +\frac{2}{\epsilon}g_K'(0) ,
 \end{align*}
and the middle summand approaches $0$ as $\epsilon\to 0$.

It remains to compute $g'_K(0)$. Note that \[g_K(-t)= \int_{S(x,\xi)=-t}\langle dx,\xi \rangle =- \int_{S(\xi,x)=t}\langle d\xi,x\rangle =-g_{K^\vee}(t).\]

Parametrizing $\alpha=\beta+s_t(\beta)$, we have
\[t^2=\langle x(\alpha), \xi(\beta)\rangle =\frac12 s_t(\beta)^2\kappa(\beta)+O(s_t^3).\]
Thus \[s_t(\beta)=\frac{\sqrt2}{\sqrt{\kappa(\beta)}}t+ O(t^2),\] and writing  
\[g_K(t)= \int_0^L \langle x'(\beta +s_t(\beta)),x(\beta)\times x'(\beta)\rangle d\beta,\]
\[\langle x'(\beta +s_t(\beta)),x(\beta)\times x'(\beta)\rangle= \sqrt{2\kappa(\beta)}t+O(t^2),   \]
we find
\[g_K'(0)=g_{K^\vee}'(0)= \sqrt 2\int_0^L\sqrt{\kappa(\beta)}d\beta. \]
\endproof

\bibliographystyle{abbrv}
\bibliography{bibliography}

\end{document}